\documentclass[11pt]{article}

\usepackage{latexsym,amsfonts}
\usepackage{amssymb,amsmath}
\usepackage{graphicx}
\usepackage[normal]{subfigure}
\usepackage{multicol}
 \usepackage{color}
\definecolor{dred}{rgb}{0.92,0,0}
\definecolor{dgreen}{rgb}{0,0.92,0}
\definecolor{dblue}{rgb}{0,0,0.92}
\definecolor{dyellow}{rgb}{0.95,0.95,0}

 \usepackage{color}
\definecolor{dred}{rgb}{0.92,0,0}
\definecolor{dgreen}{rgb}{0,0.92,0}
\definecolor{dblue}{rgb}{0,0,0.92}
\definecolor{dyellow}{rgb}{0.95,0.95,0}

\newcommand{\ds}{\displaystyle}

\def\R{\mathbb R}

\renewcommand{\div}{{\rm div}\,}

\newcommand{\curls}{curl\,}

\newcommand{\Hvec}{\mathbf{H}}

\newcommand{\Vvec}{\mathbf{V}}

\def\nvec{{\bf n}}

\def\uvec{{\bf u}}
\def\vvec{{\bf v}}
\def\xvec{{\bf x}}

\usepackage[bbgreekl]{mathbbol}

\def\Amat{\mathbb A}
\def\Mmat{\mathbb M}
\def\Tmat{\bbtau}

\newtheorem{e-proposition}[theorem]{Proposition}

\newtheorem{e-definition}[theorem]{Definition\rm}
\newtheorem{remark}{\it Remark\/}

\graphicspath{
{/Users/franckassous/Documents/Franck/Documents/Etudiants/PhD_MosheLin/Article1/}
{/Users/franckassous/Documents/Franck/Documents/Etudiants/PhD_MosheLin/Article1/Figures-article1}
{./CRAS_envoi-Moshe+code/cras_10-11-18_avec-3-objets/}
{../CRAS_Meca/CRAS_envoi-Moshe+code/cras_10-11-18_avec-3-objets/}
{./Figures/}
{./Figures-article1/}
{./CRAS_Meca/CRAS V2-resoumis/}
{/Users/franckassous/Documents/Franck/Documents/Etudiants/PhD_MosheLin/Article1/New-numerical-results-for-article1-12-01-20/draft_files_2020}
}
\title{Time Reversal for elastic scatterer location from Acoustic Recording}
\author{F. Assous\thanks{Ariel University, Ariel, Israel},
M. Lin\thanks{Ariel University, Ariel, Israel}
}
\date{}

\begin{document}
\maketitle

\bibliographystyle{plain}

\begin{abstract}
\noindent 
The aim of this paper is to study the feasibility of time-reversal methods in a non homogeneous elastic medium, from data recorded in an acoustic medium. We aim to determine, from partial aperture boundary measurements, the presence and some physical properties of elastic unknown ``inclusions", i.e. not observable solid objects, located in the elastic medium. We first derive a variational formulation of the acousto-elastic problem, from which one constructs a time-dependent finite element method to solve the forward, and then, the time reversed problem. Several criteria, derived from the reverse time migration framework, are then proposed to construct images of the inclusions, and to determine their locations. The dependence/sensitivity of the approach to several parameters (aperture, number of sources, etc.) is also investigated. In particular, it is shown that one can  differentiate between a benign and malignant close inclusions. This technique is fairly insensitive to noise in the data.
\vskip 0.5\baselineskip
%
{\small \it \noindent keywords: Time reversal;  wave propagation;  inverse problem; acousto-elastodynamics; numerical simulation}
\end{abstract}


\section*{Introduction} 
\label{Introduction}

Time reversal (TR) is a subject of very active research for over two decades, and was  experimentally developed by M. Fink in 1992 in acoustics, showing very interesting features \cite{FWCM91}. It remains an active subject of research, whether in the theoretical, physical  or numerical points of view, with many applications such as tumor detection \cite{KR2005} or earthquake prevention, see \cite{LMFCTC06} and references therein.\\

\noindent TR is a  procedure based on the reversibility property of wave propagation phenomena in non-dissipative media, like in  acoustic, electromagnetic or elastodynamics. A consequence of this property is that one can ``time-reverse" developed signals, by letting them propagate back in time to the location of the source (or scatterers) that emitted them originally. This remarkable property allows proposing innovative methods, for instance in medical imaging \cite{TOC12, S2004, B2011, CDTS2010}, seismic inversion \cite{SDMT2009} and active or passive detection problems \cite{BS2009}. The outstanding aspect of this experiment is the possibility to refocus in a very precise way the signal, without knowing the details of the source that emitted it. The first mathematical analysis can be found in \cite{Bardos02} for a homogeneous medium and in \cite{CF97}, \cite{BPZ02} for a random medium.\\

\noindent TR can also be used to identify scatterers, using the fact that a scatterer impinged upon by an incident wave behaves as a source for the reflected field. More generally, scatterers behave like secondary sources, so that TR causes reflected signals to refocus at the scatterer location.\\

\noindent TR can also be considered as an inverse problem, for which the solution is obtained without the need of an iterative process. Theoretically, solving a time reversed problem under ideal circumstances, this should yield the exact solution. However, there is always the possibility that under some (realistic!) conditions, the time reverse process will fail. This may happen due to several reasons: measurement noise, availability of only partial information in space or in time, and lack of knowledge about the medium properties or the source.\\

\noindent From a computational point of view, TR involves advancing the numerical solution of the relevant wave problem ``backward in time". The most direct use of the refocusing property allows one to solve the inverse problem of finding the location of a source when remote measurements of the wave field are given. Frequent applications are in seismology, for locating the epicenter of an earthquake from measurements taken on the ground \cite{GiTu12, LTG15}.\\

\noindent Here we are concerned with the refocusing on the scatterers of the diffracted wave to try to obtain elastic properties of an inclusion. Detection of elastic properties of a given objects inside a given ``noisy" medium, by ultrasound waves can have a lot of applications, such as  medical imaging, material detection inside a rock, see \cite{KFS2013, CDDIDRS2008} and references therein. \\

\noindent As usual for inverse problems,  this will require to define criteria, like cost functions in inverse problems \cite{GiTu12}, to measure or at least evaluate the quality of the final result, namely the refocusing of wave on the diffracting object we aim to recover.\\

\noindent The main goal of this article is to derive and analyze the TR approach for the time-dependent acousto-elastodynamics configuration. As described in \cite{B2011, DP2014} the recording of reflected waves are carried out in the fluid part of the domain, which means that only acoustic pressure is recorded. Our goal here is to deal with this partial data (acoustic pressure without knowledge of elastic displacements). The main applications of linear elastodynamics, or acousto-elastodynamics are generally structural engineering  \cite{BT2008}, where the main purpose is to examine beams that support buildings and bridges, or  seismology, where the waves can be excited naturally by earthquakes or artificially using manmade explosions, to explore the subsoil properties.  In this paper, we will focus on another application: instead of waves that excited on earth, we will consider waves that propagate through a tissue, for example, through breast tissue. The applications we have in mind are for instance the detection of tumors.\\

\noindent As far applications are concerned, this approach is supposed to determine the elastic properties of obstacles in media that are assumed to mimic breast tissue, benign and malignant tumors \cite{OCHYL1991, DP2014}. Applications we have in mind are to help the detection of  breast tumors, by using  efficient simulations of ultrasound waves. This technique can be used instead of intrusive diagnosis that performed in cases where the common diagnosis of the cancerous tissue is not guaranteed (e.g. Mammography, human palpation, Ultrasound, etc.).\\

\noindent For these reasons, mathematical models and then numerical simulations will mimic the process of sending ultrasound waves through a breast tissue according to a (partially) known mechanical properties and noise and try to detect cancerous tumors by the reflected waves, in case of partial information. Hence in this article\footnote{a short report appeared in the C. R. Mecanique (2019) \cite{AsLi19}}, in order to be closer to what happens in many real cases as in medical imaging and other fields, we do not assume that the line of receivers encloses the bounded domain $\Omega$,  we rather consider that the aperture is reduced. Then, we will apply the derived method to identify an ``inclusion" or to differentiate between a single inclusion and a two close inclusions each one with different elastic properties that correspond to different breast tumors (benign and malignant). Typically, a benign tumor corresponds to normal breast tissues, with a Young modulus between 1 and 70 KPa, whereas malignant tumors have a Young modulus varying from 15 to 500 KPa (see for instance \cite{Fern17}). Note also that,  as pointed out for instance in \cite{AKN12}, the method does not require {\em a priori} knowledge of the physical properties of the inclusions.\\

\noindent The paper is organized as follows: we present in Section \ref{Forward}, the forward problem, which is practically used only for numerically generating the measurement data (to be inverted). In Section \ref{TimeReversal}, we describe the process of time reversal (TR) and the way it is derived for the acousto-elastic problem. In Section \ref{Numerical}, we describe the cost functional used to evaluate the quality of identification. We then present numerical experiments and results, with a particular focus on breast cancer detection applications. A conclusion is drawn at the end of the paper.

\section{The forward model}
\label{Forward}

\subsection{Governing equations}

\begin{figure} 	
	\centering
	 \includegraphics[height=0.350\textwidth]{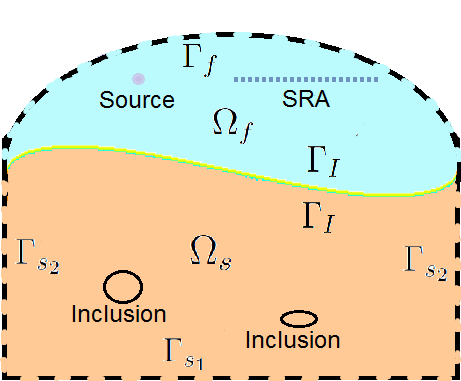}
	\caption{Geometry of a fluid-solid domain $\Omega$, with the inclusion(s) $D$ and the source-receivers array (SRA) (dashed line).}
	\label{fig:simulation_domain}
\end{figure}

\noindent We consider a non homogeneous two-dimensional fluid-solid domain $\Omega\subset  \R^N, N=2$, that consists of a top-most acoustic layer 
$\Omega_f$ and, underlying it, an elastic domain $\Omega_s$. The acoustic part $\Omega_f$ corresponds to a homogeneous fluid, characterized by a known wave propagation velocity $V_f^p=\sqrt{\frac{\lambda_f}{\rho_f}}$,  $\rho_f$ and $\lambda_f$ denoting respectively the density and the bulk modulus of the fluid. In order to be closer to realistic cases as in geophysics or medical imaging, we consider that the elastic part of the domain $\Omega_s$ corresponds to a layered solid, characterized for each layer by a known wave propagation velocities $V_s^p=\sqrt{\frac{\lambda_s+2\mu_s}{\rho_s}}$ and $V^s_s=\sqrt{\frac{\mu_s}{\rho_s}}$, $\rho_s$ and $\lambda_s, \mu_s$ being respectively the density and the Lam\'e coefficients of the considered layer (see Figure \ref{fig:simulation_domain}).\\

\noindent Hence, the incident wave, generated by a known source, or being an  incoming wave from infinity, is first propagated in the fluid part of  the medium,  then transmitted (and partially reflected) into the solid medium. This wave is then scattered  by a buried obstacle D located inside the solid domain (which may contain one or several inclusions) characterized by different physical properties from the surrounding elastic layer. The resulting scattered wave is then recorded by a source-receivers array (SRA) located in the fluid part of the domain (see Figure \ref{fig:simulation_domain}). We also assume that after a time $T_f$, the total field is negligible in the domain of interest $\Omega$.\\

\noindent To model this phenomenon, we introduce the scalar wave equation in the fluid part $\Omega_f \subset \Omega$, and the elastodynamic wave equation in the solid part $\Omega_s \subset \Omega$. The fluid-solid interface $\Gamma_I$ (see Figure \ref{fig:simulation_domain}) will be handled by interface transmission conditions. 

\begin{remark}
\noindent We present here the coupled system of acousto-elastodynamics equations. The (fully) elastodynamics equations can be derived, by simply removing the fluid part (and related equations) from the model described below, and by adding a term source (and recorders) in the elastic domain. Note also that  there is no "conceptual" difficulty to derive the same problems in three-dimensional geometry.
\end{remark}

\noindent Let us consider first the acoustic problem. We denote by $\partial \Omega_f$ the boundary of $\Omega_f$ and $\textbf{n}=(n_1,n_2)$ is the the unit outward normal vector to $\partial \Omega$. We denote by $p(\xvec,t)$ the pressure  on a time $t$, $\xvec=(x_1,x_2) \in \Omega_f$, and by $f(\xvec,t)$ a given source. In these conditions,  the acoustic wave equation in $\Omega_f$ can be written as
\begin{equation} \label{acoustic_wave1}
\frac{{1}}{{\lambda_f}}\frac{{\partial ^2 p}}{{\partial t^2}} - 
\div(\frac{1}{\rho_f}\nabla p)
= f\,,
\end{equation}
together with homogeneous initial conditions at the initial time $t=0$
\begin{equation} \label{acoustic_wave2}
p(t=0) = 0 , \quad \frac{{\partial p}}{{\partial t}}(t=0) = 0\,.
\end{equation}
Let us consider now that the boundary $\partial \Omega_f$ can be split in $\partial \Omega_f=\Gamma_f \cup \Gamma_I$, 
$\Gamma_I$ being the fluid-solid interface introduced above, and $\Gamma_f$ the artificial outer surface of $\Omega_f$. Assuming that the boundary $\Gamma_f $ is a ball, we supplement Eqs. (\ref{acoustic_wave1}-\ref{acoustic_wave2}) with the first order Bayliss-Turkel absorbing boundary conditions BT$^1$ \cite{BaTu80,BaGT82}
\begin{equation} \label{acoustic_wave3BT}
\frac{\partial p}{\partial t} + V_p \frac{\partial p}{\partial r} + V_p \frac{p}{2r} = 0 \quad \mbox{ on } \Gamma_f\,,
\end{equation}
$r$ denoting the radial coordinate. 
\begin{remark}
\noindent In the case of straight boundary $\Gamma_f $, we will use the Engquist-Majda absorbing boundary condition \cite{ClEn77} instead of (\ref{acoustic_wave3BT})
\begin{equation} \label{acoustic_wave3}
\frac{{\partial p}}{{\partial t}} = -V_p\, 
\nabla p\cdot \textbf{n} \quad \mbox{ on } \Gamma_f\,.
\end{equation}
\end{remark}
On the interface $\Gamma_I$, transmission conditions will be added, see (\ref{acoustic_continuity}) below.\\

\noindent Similarly for the solid part $\Omega_s$ of the domain, the velocity $\textbf{u}(\xvec,t)=(u_1(x_1,x_2,t), u_2(x_1,x_2,t))$ at a point $\xvec=(x_1,x_2) \in \Omega_s$ satisfies, for all $t>0$,
\begin{equation} \label{elastodynamic_wave3}
\rho_s\frac{\partial ^2 u_i }{\partial t^2} - \sum_{j=1,2} \frac{\partial}{\partial x_j} \tau_{ij}(\uvec)= 0\,,\quad i=1,2\,,
\end{equation}
where $\Tmat(\uvec):=(\tau_{i,j}(\uvec))_{1 \leq i,j \leq 2}$ is the (time derivative of the) classical stress matrix defined by 
$$
\tau_{ij}(\uvec) = \lambda_s \div\uvec \, \delta_{ij}+ \mu_s(\ds\frac{{\partial u_i}}{{\partial x_j}} + \frac{{\partial u_j}}{{\partial x_i}}), \quad
i,j = 1,2 \,.
$$
\begin{remark}
These equations can be equivalently formulated in a vector form as  
$$
\rho_s\frac{{\partial ^2 \textbf{u}}}{{\partial t^2}} - \div\Tmat = 0\,, \quad \mbox{ or } \quad
\rho_s\frac{{\partial ^2 \textbf{u}}}{{\partial t^2}} - \div(\mu_s\nabla \textbf{u}) -\nabla((\lambda_s + \mu_s)\nabla\cdot \textbf{u})= 0\,.
$$
\end{remark}

\noindent As for the fluid part of the domain, this equation is supplemented with homogeneous initial conditions
\begin{equation} \label{elastodynamic_wave4}
\textbf{u}(t=0) = 0, \quad \frac{{\partial \textbf{u}}}{{\partial t}}(t=0) = 0\,.
\end{equation}
Here also, the boundary of the solid domain $\partial \Omega_s$ can be split in $\partial \Omega_s=\Gamma_s \cup \Gamma_I$, where $\Gamma_s$ denotes the artificial outer surface of $\Omega_s$. Equations (\ref{elastodynamic_wave3}) and (\ref{elastodynamic_wave4}) are supplemented with absorbing boundary conditions on $\Gamma_s$, as proposed for instance in \cite{Halp80}, \cite{BJRT85}.\\

\noindent  For our purpose here, let us introduce the first order absorbing condition that can be written here as 
\begin{equation}\label{ABCElast_general}
\Tmat(\uvec)\, \nvec + \Mmat \frac{\partial \uvec}{\partial t}=0  \quad \mbox{ on } \Gamma_s\,,
\end{equation}
where $\Mmat$ is a $N\times N$ matrix  (here $N=2$) defined by
$$
\Mmat =
\left(\begin{array}{cc}
n_1 & n_2\\
n_2 & -n_1
\end{array}
\right)
\left(\begin{array}{cc}
\rho V_p & 0\\
0 & \rho V_s
\end{array}
\right)
\left(\begin{array}{cc}
n_1 & n_2\\
n_2 & -n_1
\end{array}
\right)\,,
$$
that is, a symmetric and positive definite matrix. Physically, this condition expresses that the normal stress is inversely proportional to the velocity. In particular, one can define a total energy as a non increasing function of the time, and condition (\ref{ABCElast_general})  is absorbing, leading to a well posed mathematical problem. A particular and very common case of Eq. (\ref{ABCElast_general}) corresponds to horizontal and vertical boundary, says $\Gamma_{s_1}$ and $\Gamma_{s_2}$ respectively. In that case, the condition is expressed as
\begin{equation} 
\label{elastodynamic_wave5}
\Amat  \frac{\partial \uvec}{\partial t} = \Tmat(\uvec)\,\nvec\,,
\end{equation}
where the matrix $\Amat$ is a diagonal $N \times N$ matrix, with $A_{11}=-\sqrt{\rho_s(\lambda_s + 2\mu_s)}, A_{22}=-\sqrt{\rho_s \mu_s}$ for horizontal boundaries, and the contrary for vertical boundaries.\\ 

\noindent Finally, we have to deal with the transmission conditions to impose on the interface $\Gamma_I$ between the fluid and the solid parts. 
\begin{enumerate}
\item We begin with the continuity of the normal (with respect to $\Gamma_I$) component,  that reflects that the jump across $\Gamma_I$ of the displacement is equal to 0, and so is the jump of the velocity $\uvec$, that is $[\uvec\cdot \nvec]_{\Gamma_I}=0$. Now, using the Euler equation written in velocity-pressure $\rho_f\ds\frac{\partial \uvec}{\partial t}+\nabla p=0$,  and differentiating in time the jump relation across $\Gamma_I$, we readily obtain the first relation
\begin{equation} \label{acoustic_continuity}
\frac{1}{\rho_f}  \frac{\partial p}{\partial \nvec} =\frac{1}{\rho_f}\nabla p\cdot \nvec =- \frac{\partial \uvec}{\partial t}\cdot \nvec\,, 
\end{equation} 
$\textbf{n}$ denoting here the unit outward normal vector from $\partial \Omega_s$ to $\partial \Omega_f$.
\item The second relation to impose is  the continuity of the normal component of the stress tensor, that gives, agaim after differentiation in time
\begin{equation} \label{elastic_continuity1}
 \frac{{\partial p}}{{\partial t}} n_i= \sum_{j=1,2}\tau_{ij}(\uvec) \,n_j, \quad i=1,2\,, 
\end{equation} 
or 
$
\ds \frac{{\partial p}}{{\partial t}} \nvec= \Tmat(\uvec) \,\nvec 
 $
in a vector form.
\end{enumerate}

\begin{remark}\label{transhoriz}
\noindent For an horizontal interface, relation (\ref{acoustic_continuity}) can be simply expressed 
$
\ds\frac{1}{\rho_f} \frac{{\partial p}}{{\partial x_2}} = -\ds\frac{{\partial u_2}}{{\partial t}}\, ,
$
whereas  relation (\ref{elastic_continuity1}) is written as 
$
0=\ds\frac{{\partial u_1}}{{\partial x_2}} + \frac{{\partial u_2}}{{\partial x_1}}, \quad 
\ds\frac{{\partial p}}{{\partial t}} = \lambda_s \frac{{\partial u_1}}{{\partial x_1}} + (\lambda_s + 2\mu_s)\frac{{\partial u_2}}{{\partial x_2}}.
$
\end{remark}

\noindent As a consequence of the transmission conditions, there is conservation of $N+1=3$ scalar quantities. Note also that these conditions appear as {\em natural boundary conditions} in the variational formulation that will be derived below, due to the use of the pressure-velocity formulation (\ref{acoustic_wave1}-\ref{elastodynamic_wave3}).\\

\subsection{Variational formulation}
\label{VF-Forward}

\noindent We now introduce the variational formulation for the acousto-elastodynamics equations, whose derivation is not  classical, due to the coupling between the fluid and the solid parts. This is why the elastodynamics wave equation is expressed in velocity (and not in displacement). This approach leads to an "easy" coupling on the interface $\Gamma_I$, and allows to perform the computations in a ÒnaturalÓ and expected 
way.\\

\noindent Let us take a test function $q \in V$ for acoustic wave equation,  $V$ being an {\em ad hoc} Sobolev space. Classically,  $V=H^1(\Omega_f):=\{v \in L^2(\Omega_f);  \nabla v \in L^2(\Omega_f)\}$. Multiplying (\ref{acoustic_wave1}) by $q\in V$, integrating over the domain $\Omega_f$, and using a classical integration-by-parts formula, we get
\begin{eqnarray*} 
&&\int_{\Omega_f}\frac{{1}}{{\lambda_f}}\frac{{\partial ^2 p}}{{\partial t^2}}q \,d\omega 
+\int_{\Omega_f}\frac{1}{\rho_f}\nabla p \cdot \nabla q \,d\omega
- \int_{\partial \Omega_f}\frac{1}{\rho_f} \frac{\partial p}{\partial \nvec} q \,d\sigma 
= \int_{\Omega_f}fq \,d\omega\,,
\end{eqnarray*}
\noindent Taking into account that  $\partial \Omega_f=\Gamma_f \cup \Gamma_I$, one uses the continuity of the normal component (\ref{acoustic_continuity}) on $\Gamma_I$, and the ABC (\ref{acoustic_wave3BT}) on $\Gamma_f$ (or the ABC (\ref{acoustic_wave3}) the case occurring). We get
\begin{eqnarray} 
\label{acoustic_vf_v1}
&&\hspace*{-1.8cm}\int_{\Omega_f}\frac{{1}}{{\lambda_f}}\frac{{\partial ^2 p}}{{\partial t^2}}q \,d\omega 
+\int_{\Omega_f}\frac{1}{\rho_f}\nabla p \cdot \nabla q \,d\omega
+ \int_{\Gamma_I} \frac{\partial\uvec}{\partial t}\cdot \nvec \, q \,d\sigma  \nonumber \\
&&\hspace*{1.cm}+ \int_{\Gamma_f} \left(\frac{1}{\sqrt{\rho_f\lambda_f}}\frac{{\partial p}}{{\partial t}}+\frac{1}{\rho_f} \frac{p}{2r}\right)
 \,q \,d\sigma 
= \int_{\Omega_f}f q \,d\omega\,,
\end{eqnarray} 
where the second term in the integral on $\Gamma_f$, that is $\ds\int_{\Gamma_f} \frac{1}{\rho_f} \frac{p}{2r} \,q \,d\sigma$, is replaced by 0 if one uses the ABC (\ref{acoustic_wave3}).\\

\noindent Similarly for elastodynamics wave equation, to derive the variational formulation in $\Omega_s$, we consider $\vvec=(v_1,v_2)$ a vector test function, that belongs to the relevant Sobolev space $\Vvec:=\Hvec^1(\Omega_s)$. We use the Green's identity (see \cite{Ciarlet1}, p. 288) 
\begin{equation}
\label{eq:GreensIdentity1}
-\int_{\Omega_s} \div\mathbb{S}\cdot\rm{\textbf{v}}\,d\omega=\int_{\Omega_s} \mathbb{S}:\nabla\rm{\textbf{v}}\,d\omega-\int_{\partial\Omega_s} \left(\mathbb{S}\,\,\rm{\textbf{n}}\right)\cdot\rm{\textbf{v}}\,d\sigma\,, \quad \forall \vvec \in \Vvec\,,
\end{equation}
where $\mathbb{S}$ denotes a symmetric tensor field and the symbol $:$ stands for the contracted product of two tensors namely, $\bbsigma(\uvec): \bbsigma(\vvec)=\ds\sum_{i=1}^{2}\sum_{j=1}^{2} \bbsigma_{ij}(\uvec)\bbsigma_{ij}(\vvec)
= \bbsigma_{ij}(\uvec)\bbsigma_{ij}(\vvec)$ with the Einstein summation convention.\\

\noindent Multiplying (\ref{elastodynamic_wave3}) by $\vvec=(v_1, v_2) \in \Vvec$, integrating over the domain $\Omega_s$, and using integration by parts (\ref{eq:GreensIdentity1}),  we get
\begin{eqnarray*} 
&&\hspace*{-1.2cm}\int_{\Omega_s} \rho_s \frac{\partial ^2 \uvec}{\partial t^2} \cdot \vvec \,d\omega 
+\int_{\Omega_s}\lambda_s \div\uvec\, \div\vvec \,\,d\omega\,
+\int_{\Omega_s} 2 \mu_s \Tmat(\uvec) : \Tmat(\vvec) \,d\omega
-\int_{\partial \Omega_s} \tau_{ij}(\uvec) \,v_i\, n_j \,d\sigma=0\,.
\end{eqnarray*} 

\noindent Using now the continuity of the normal component of the stress tensor (\ref{elastic_continuity1}) on the interface, and the ABC (\ref{ABCElast_general}) on $\Gamma_s$ (or the ABC (\ref{elastodynamic_wave5}) the case occurring), we obtain
$$
\int_{\Omega_s} \rho_s \frac{\partial ^2 \uvec}{\partial t^2}\cdot \vvec \,d\omega 
+\int_{\Omega_s}\lambda_s \div\uvec\, \div\vvec \,
+ 2 \mu_s \tau_{ij}(\uvec) \, \tau_{ij}(\vvec) \,d\omega
-\int_{\Gamma_I} \frac{{\partial p}}{{\partial t}} \,\vvec\cdot \nvec \,d\sigma\, 
+\int_{\Gamma_{s}}\Mmat  \frac{\partial \uvec}{\partial t} \cdot\vvec \,d\sigma=0\,.
$$
where the integral $\ds\int_{\Gamma_{s}}\Mmat  \frac{\partial \uvec}{\partial t} \cdot\vvec \,d\sigma$ is replaced by $\ds-\int_{\Gamma_{s}}\Amat  \frac{\partial \uvec}{\partial t} \cdot\vvec \,d\sigma$
when using the ABC (\ref{elastodynamic_wave5}) instead of (\ref{ABCElast_general}).\\

\begin{remark}
\noindent For an horizontal interface parallel to the $x_1$ axis, the variational formulation (\ref{acoustic_vf_v1}) can be simply expressed 
\begin{eqnarray*} 
&&\hspace*{-1.cm}\int_{\Omega_f}\frac{{1}}{{\lambda_f}}\frac{{\partial ^2 p}}{{\partial t^2}}q \,d\omega 
+\int_{\Omega_f}\frac{1}{\rho_f}\nabla p \cdot \nabla q \,d\omega
 +
\int_{\Gamma_{I}}\frac{{\partial u_2}}{{\partial t}} \, q \,d\sigma  
+ \int_{\Gamma_f}\left(\frac{1}{\sqrt{\rho_f\lambda_f}}\frac{\partial p}{\partial t} +\frac{1}{\rho_f} \frac{p}{2r}\right)\, q \,d\sigma = \int_{\Omega_f}f q \,d\omega\,.
\end{eqnarray*} 
In addition, in the common case of horizontal and vertical absorbing boundaries $\Gamma_{s_1}$ and $\Gamma_{s_2}$, one has 
\begin{eqnarray*} 
&&\hspace*{-2.5cm}\int_{\Omega_s}\rho_s\frac{\partial ^2 u_1}{\partial t^2}v_1 \,d\omega 
 + \int_{\Omega_s} ((\lambda_s + 2\mu_s)\frac{\partial u_1}{\partial x_1} + \lambda_s \frac{\partial u_2}{\partial x_2})\frac{\partial v_1}{\partial x_1} 
 + \mu_s (\frac{\partial u_1}{\partial x_2} + \frac{\partial u_2}{\partial x_1})\frac{\partial v_1}{\partial x_2} \,d\omega  \\
&&\hspace*{1.5cm}+ \int_{\Gamma_{s_{1}}}\sqrt{\rho_s\mu_s}\frac{{\partial u_1}}{{\partial t}}v_1 \,d\sigma
+ \int_{\Gamma_{s_{2}}}\sqrt{\rho_s(\lambda_s + 2\mu_s)}\frac{{\partial u_1}}{{\partial t}}v_1 \,d\sigma=0 
\end{eqnarray*} 
\begin{eqnarray*} 
&&\hspace*{-0.6cm}\int_{\Omega_s}\rho_s\frac{{\partial ^2 u_2}}{{\partial t^2}}v_2  \,d\omega 
+ \int_{\Omega_s}  ((\lambda_s + 2\mu_s)\frac{\partial u_2}{\partial x_2} + \lambda_s \frac{\partial u_1}{\partial x_1})\frac{\partial v_2}{\partial x_2} 
+ \mu_s (\frac{\partial u_1}{\partial x_2} + \frac{\partial u_2}{\partial x_1}) \frac{\partial v_2}{\partial x_1}\,d\omega \\
&&\hspace*{2.9cm}+ \int_{\Gamma_{s_{2}}}\sqrt{\rho_s\mu_s}\frac{{\partial u_2}}{{\partial t}}v_2 \,d\sigma
+ \int_{\Gamma_{s_{1}}}\sqrt{\rho_s(\lambda_s + 2\mu_s)}\frac{{\partial u_2}}{{\partial t}}v_2 \,d\sigma 
   - \int_{\Gamma_{I}}\frac{{\partial p}}{{\partial t}} v_2\,d\sigma = 0
\end{eqnarray*} 
\end{remark}

\begin{remark}
The transmission conditions at the interface $\Gamma_I$ appear through an anti-symmetric form $c(\cdot,\cdot)$ defined by
$$
\frac{d}{dt} c((p, \uvec),(q,\vvec)):=\frac{d}{dt} \left(\int_{\Gamma_I} - p \,\vvec\cdot \nvec \,d\sigma\,
+ \int_{\Gamma_I} \uvec\cdot \nvec \, q \,d\sigma \right)
$$ 
This allows us to define an energy of the system and to prove that this energy is conserved (in the absence of ABC), or decreases (in the presence of ABC). Let us give here the sketch of the proof:\\

\noindent We first introduce the following notations:
\begin{eqnarray*}
&&((p, \uvec),(q,\vvec))=\int_{\Omega_f}\frac{1}{\lambda_f} p \, q \,d\omega
+ \int_{\Omega_s} \rho_s \uvec \cdot \vvec \,d\omega\,,\\
&&a((p, \uvec),(q,\vvec))=a_f(p,q)+a_s(\uvec,\vvec)=
\int_{\Omega_f}\frac{1}{\rho_f}\nabla p \cdot \nabla q \,d\omega
+\int_{\Omega_s}\lambda_s \div\uvec\, \div\vvec \,+ 2 \mu_s \tau_{ij}(\uvec) \, \tau_{ij}(\vvec) \,d\omega\,,\\
&&b((p, \uvec),(q,\vvec))= 
\int_{\Gamma_{f}} \frac{1}{\sqrt{\rho_f\lambda_f}} p \,q \,d\sigma 
+\int_{\Gamma_{s}}\Mmat  \uvec \cdot\vvec \,d\sigma\,,\mbox{ and } e(p,q) = \int_{\Gamma_{f}} \frac{1}{\rho_f} \frac{p}{2r} \,q \,d\sigma\,, 
\end{eqnarray*}
so that the variational formulation can be written as, after extinction of the source $f(\xvec,t)$
$$
\frac{d^2}{dt^2}((p, \uvec),(q,\vvec))+a((p, \uvec),(q,\vvec))+\frac{d}{dt}b((p, \uvec),(q,\vvec))+\frac{d}{dt}c((p, \uvec),(q,\vvec))+ e(p,q)=0\,.
$$

\noindent To define the energy of the system, we choose $q=\ds\frac{\partial p}{\partial t}$ and $\vvec=\ds\frac{\partial \uvec}{\partial t}$ in the formulation above. This yields the antisymmetric term $c$ to vanish, and we get
\begin{equation}
\label{energy}
\frac{1}{2}\frac{d}{dt}((\frac{\partial p}{\partial t}, \frac{\partial \uvec}{\partial t}),(\frac{\partial p}{\partial t}, \frac{\partial \uvec}{\partial t}))+\frac{1}{2} a((p, \uvec),(p,\uvec))= 
-b((\frac{\partial p}{\partial t}, \frac{\partial \uvec}{\partial t}),(\frac{\partial p}{\partial t}, \frac{\partial \uvec}{\partial t}))
-\frac{1}{2}  e(p,p)\,.
\end{equation}

\noindent The first term of (\ref{energy}) can be interpreted as a kinetic energy, whereas the second one is a term of deformation. Without term of ABC, namely for $b \equiv e \equiv 0$, the energy defined above is conserved, whereas it decreases in presence of ABC.
\end{remark}

\section{Recreate the past with time reversal } 
\label{TimeReversal}

\noindent We wish to recover an unknown obstacle buried in the elastic part of the medium from field measurements, by using the TR method. Basically, it takes advantage of the reversibility of wave propagation phenomena in a non dissipative unknown medium to back-propagate signals to the sources that emitted them \cite{AKNT11}. In this section, we will first present the principle of the Time Reversal techniques, then derive it for the fluid-solid problems.

\subsection{Principle of time reversal}
\label{TRprinciple}
Let us describe briefly the principle of TR, for a general propagation equation. We consider an incident wave $u^I$ impinging an obstacle $D$ characterized by different physical properties from the surrounding medium. 
Assuming that the total field $u$ satisfies a linear hyperbolic equation or system of equations, says $\mathcal{L}$, one can write
\begin{equation} \label{tr_principle_initial}
\begin{gathered}
\mathcal{L}(u) = 0 \: \mbox{ in } \: R^N\,,
\end{gathered}
\end{equation}
together with zero initial conditions.\\

\noindent Now, we introduce a bounded domain $\Omega$ that surrounds $D$, and assume that the incident wave $u^I$ is generated by a source such that after a time $T_f$ the total field $u$ is negligible in $\Omega$. The time-reversed solution $u^R$ of equation (\ref{tr_principle_initial}) is simply defined as $u^R:=u(T_f-t)$. Denoting by $\Gamma_{SRA}$ the boundary of source-receivers array on which the scattered field is recorded.  We impose Dirichlet boundary conditions on $\Gamma_{SRA}$ equal to the time reversal of the recorded fields and zero initial conditions.\\

\noindent  In real configurations, since the physical properties of the inclusion or its exact location are unknowns, the form of the operator $\mathcal{L}$ is unknown inside the inclusion $D$, but only outside $D$. In these conditions, the corresponding $u^R$ satisfies equation (\ref{tr_principle_initial}) (reversed in time) only outside the object $D$, namely in $\Omega \backslash D$. Hence, one can not define a well-posed boundary value problem on $u^R$ in $\Omega \backslash D$.\\

\noindent Nevertheless, the classical TR approach overcomes this difficulty by solving the problem in the entire domain $\Omega$, assuming that there is no inclusion $D$, building in this way an ``approximate" time-reversed solution. More precisely, denoting by $w^R$ this ``approximate" time-reversed solution, one solves in the entire domain $\Omega$, i.e. without inclusion $D$,
$$
\mathcal{L}_0(w^R) = 0 \: \mbox{ in } \Omega\,,
$$
together with Dirichlet boundary conditions on $\Gamma_{SRA}$, equal to the time reversal of the recorded fields, and zero initial conditions for the TR problem.

\subsection{Time Reversal for acousto-elastodynamics equations}
\label{TRacoustoElastic}

\noindent Considering a fluid-solid problem as described above, we aim to reconstruct the time reversed acousto-elastic wave $(p^R,\uvec^R)$, scattered from a scatterer $D$, which may contain several obstacles or inhomogeneities. As explained above, our measurement data have been recorded on the boundary $\Gamma_{SRA}$, from a known {\em incident} wave.\\

\noindent In the configurations considered in the present paper, due to the non-homogeneity of the fluid-solid domain, made of two or several layers, the {\em incident} field is defined as the pressure-velocity fields that propagate in the absence of scatterer $D$. This assumes that the location of the interface $\Gamma_I$ and of the eventual other discontinuities is known. However, the recording being performed in the fluid part, we have likewise to solve a time reversed fluid-solid problem only from acoustic recording, that is from the measurement of the pressure $p$ in the fluid part.\\

\noindent We describe now the TR system of equations, together with the boundary and interface conditions. The resulting system being quite similar to the one of the forward problem, we will explain briefly the construction of the variational formulations.\\

\noindent Consider first $\Omega_f$ the fluid part of the domain, where the time-reversed pressure $p^R(\xvec,t')$ is defined by $p^R(\xvec,t')=p(\xvec,T_f-t)$,  $\xvec \in \Omega_f$. The wave equation involving only second order time derivatives, this definition ensures that the reverse field $p^R(\xvec,t')$ is a solution to the acoustic wave equation 
$$
\frac{{1}}{{\lambda_f}}\frac{{\partial ^2 p^R}}{{\partial t^{'2}}} - \div(\frac{{1}}{{\rho_f}} \nabla p^R)
= 0 \,,
$$
together with TR initial conditions and TR absorbing boundary conditions on $\Gamma_f$, analogous to (\ref{acoustic_wave2}) and (\ref{acoustic_wave3BT}) or (\ref{acoustic_wave3}). In that case, there is no external source $f$. Nonetheless, the boundary $\Gamma_{SRA}$ modeling a source-receivers array, a Dirichlet boundary condition equal to the time reversal of the recorded fields is imposed there, namely $p^R(t')=p(T_f-t)$ on  $\Gamma_{SRA}$.\\

\noindent Secondly, for the elastic part of the domain,  let us denote by $\textbf{u}^R(\xvec,t')=(u_1^R(x_1,x_2,t'), u_2^R(x_1,x_2,t'))$ the time-reversed velocity solution to equation (\ref{elastodynamic_wave3}), that solves
$$
\rho_s\frac{\partial ^2 u^R_i }{\partial t^{'2}} - \sum_{j=1,2} \frac{\partial}{\partial x_j} \tau_{ij}(\uvec^R)= 0\,,\quad i=1,2\,,
$$
together with (TR) initial conditions, and (TR) absorbing boundary conditions on $\Gamma_s$, with expressions analogous to  (\ref{elastodynamic_wave4}) and (\ref{ABCElast_general}) or (\ref{elastodynamic_wave5}).\\

\noindent Finally, we derive the time-reversed continuity transmission conditions at the interface $\Gamma_I$. We readily arrive at
\begin{equation}
 \label{acoustic_continuity_r1}
\frac{{1}}{{\rho_f}} \frac{{\partial p^R}}{{\partial \nvec}} =  \frac{\partial \uvec^R}{\partial t'}\cdot \nvec\,,
\end{equation}
\begin{equation} \label{tr_elastic_continuity1_general}
-\frac{{\partial p^R}}{{\partial t'}} n_i= \sum_{j=1,2}\tau_{ij}(\uvec^R) \,n_j, \quad i=1,2\,,  
\end{equation} 
Here again, one can obtain simplified versions of (\ref{acoustic_continuity_r1}) and (\ref{tr_elastic_continuity1_general}) for an horizontal interface, in the same spirit as in remark \ref{transhoriz}.\\

\noindent Let us conclude this section with the variational formulation for the time-reversed acousto-elastic equations. It is obtained by following the same principles as for the forward one. For completeness, we present below the final expression, that is written
$$
\int_{\Omega_f}\frac{{1}}{{\lambda_f}}\frac{{\partial ^2 p^R}}{{\partial t^{'2}}}q \,d\omega 
+\int_{\Omega_f}\frac{1}{\rho_f}\nabla p^R \cdot \nabla q \,d\omega
-  \int_{\Gamma_I} \frac{\partial\uvec^R}{\partial t'}\cdot \nvec \, q \,d\sigma  
- \int_{\Gamma_{f}}\frac{{1}}{{\sqrt{\rho_f\lambda_f}}}\frac{{\partial p^R}}{{\partial t'}} \,q \,d\sigma 
=  \int_{\Gamma_{SRA}}f_{SRA}\,q \,d\omega\,,
$$

$$
\int_{\Omega_s} \rho_s \frac{\partial ^2 \uvec^R}{\partial t^{'2}}\cdot \vvec d\omega 
+\int_{\Omega_s}\lambda_s \div\uvec^R \div\vvec 
+ 2 \mu_s \tau_{ij}(\uvec^R)  \tau_{ij}(\vvec) d\omega
\\ +\int_{\Gamma_I} \frac{{\partial p^R}}{{\partial t'}} \vvec\cdot \nvec d\sigma 
+\int_{\Gamma_{s}}\Amat  \frac{\partial \uvec^R}{\partial t'} \cdot\vvec d\sigma=0\,.
$$
\noindent  The forward and reversed formulations are discretized by similar numerical schemes. Both are approximated by using a $\mathrm{P}_2$ finite element method with FreeFem++ \cite{Hech12}. The time derivative is approximated by a second order centered finite difference scheme, which is time reversible also on the numerical level.

\section{Numerical results} 
\label{Numerical}

\noindent In the following, we will illustrate our approach through numerical examples,  in the case of two scatterers located in the elastic part of the medium. Let us first describe the principle of the numerical process: 
\begin{enumerate}
\item To first create the synthetic data, an incident wave is generated by a point source located in the fluid part $\Omega_{f}$,  and is propagated in the medium during a time $T_f$ such that,  for $t  \geq T_f$, the total field is almost negligible. On the boundary $\Gamma_{SRA}$ - also located in the fluid part -  the forward signal is recorded. The incident field, as defined in subsection \ref{TRacoustoElastic}, is also computed in this first step, to enable construction of the scattered field $\uvec^S$ from the total one by subtraction.
\item The time-reversed computation consists in back propagating the recorded scattered data - the difference
between the recorded total field and the incident field, both reversed in time - from the SRA in the entire domain, 
without knowing the location and the properties of the inclusions which diffracted the signals. In other words, the recorded data are back propagated in the medium without knowing the presence of the inclusions in the elastic medium.
\item Now, to determine the locations and some properties of these inclusions, one introduces a criterion that consists in the cross-correlating between the forward incident wave $\uvec^I$ with the reversed wave $\uvec^S_R$, namely
\begin{equation} \label{RTM_equation}
RTM(\xvec) = \int_{0}^{T_f} \uvec^S_R(T_f - t, \xvec) \cdot \uvec^I(t, \xvec) dt\,.
\end{equation} 
Criterion (\ref {RTM_equation}) is derived in the same spirit as those involved in reverse-time migration, in classical applications of earth imaging and images the discontinuities of the propagation speed. There are many references on the subject, especially in the framework of geophysics modeling and inversion (see, among others,  \cite{Berk84} or \cite{Clae85} and references therein), where the adjoint-solution formulation is a particular TR problem, cf. \cite{AKN12}. This method is also related to those arising from topological derivative-based methods, as described in \cite{DoGE05,BJRG15,Bonn06}. This also allows us to construct a method imaging the unknown scatterers in the medium. Note that analogous criteria can be defined using only one of the components of $\uvec$, for example for $u_2$
\begin{equation} \label{RTM_equationUY}
RTM_{{u_{2}}}(\xvec) = \int_{0}^{T_f} u_{2,R}^S(T_f - t, \xvec) \cdot u_{2}^I(t, \xvec) dt\,,
\end{equation} 
or with the $\div$ and $\curls$ operators. For instance for the divergence, one will also  consider the following imaging function
\begin{equation} 
\label{RTM_equationDIV}
RTM_{div}(\xvec) = \int_{0}^{T_f} \div\uvec^S_R(T_f - t, \xvec)\cdot \div \uvec^I(t, \xvec) dt\,.
\end{equation}
In the same way, one could also consider mathematical energies for the criterion, that can are nothing but the $L^2$ (without including the gradients) or $H^1$ (with the gradients) norms.
\end{enumerate}

\subsection{Imaging with TR: numerical results}
\label{NumResImagTR}
\begin{figure}
	\centering
	\includegraphics[width=80mm,height=48mm]{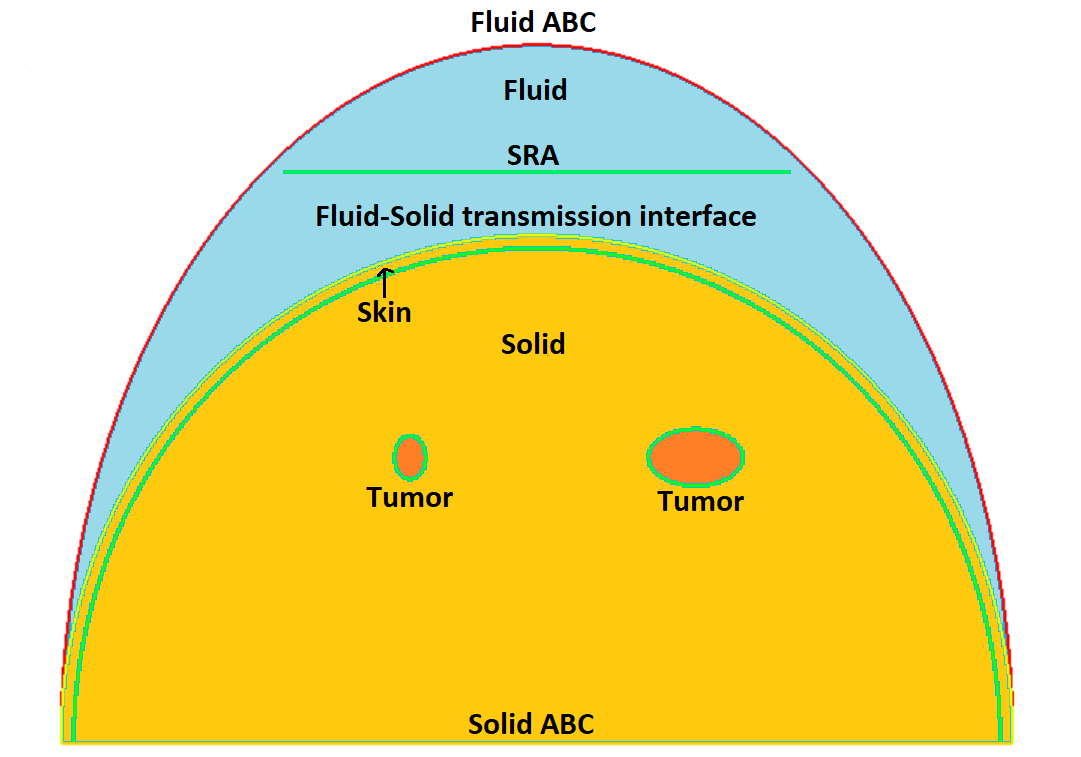}
	\caption{Geometry of typical domain $\Omega$ sketching a breast cancer detection problem, based on details taken from \cite{DBBAP2010}.}
	\label{fig:real_world_example}
\end{figure}

\noindent In the following, we will illustrate our approach through several numerical examples. Let us begin by describing the geometry of the domain. We consider the medium sketched in Figure \ref{fig:real_world_example} (inspired from \cite{DBBAP2010}), made of a fluid part (top) and of an elastic one (bottom) sketching a breast tissue geometry. It is a heterogeneous medium, as it contains a very thin layer, modelling the presence of the skin. The SRA is an horizontal line as depicted on this figure.\\

\noindent As recalled above, the first step consists in creating the synthetic data by solving the forward problem and by recording the values of the pressure $p(\xvec,t)$ on $\Gamma_{SRA}$. For this purpose, an incident wave is generated by a point source located in the fluid part, in which the density $\rho = 1000 kg/m^3$ and $\lambda = 2.25 GPa$, that is a velocity $V_p = 1500 m/sec$. For the solid part, in the very thin layer (says the skin), we take $\rho = 1150 kg/m^3 $, $\lambda = 6.66 GPa$ and $\mu = 66.66 kPa$, which correspond to compressional and shear wave velocities  $V_p = 2407.73 m/sec$ and $V_s = 7.61 m/sec$, whereas $\rho = 1000 kg/m^3$, $\lambda = 1.83 GPa$ and  $\mu = 18.33 kPa$ - corresponding to velocities  $V_p = 1354.02 m/sec$ and $V_s = 4.28 m/sec$ - in the rest of the elastic domain.\\ 
 
\noindent In addition, there are two penetrable obstacles with different size, shape, and elastic properties, located in the elastic part, outside the skin. The first one, that represents a benign tumor, is characterised by $\rho = 1000 kg/m^3 $, $\lambda = 2.16 GPa$,  $\mu = 21.66 kPa$, which correspond to P-wave and S-wave velocities $V_p = 1471.97 m/sec$ and $V_s = 4.65 m/sec$. The second one models a malignant tumor, with the same $\rho$,  $\lambda = 2.99 GPa$ and  $\mu = 30 kPa$. The corresponding P-wave and S-wave velocities are $V_p = 1732.06 m/sec$ and $V_s = 5.47 m/sec$.  The size of the scatterers is  between $0.4 \lambda_W$ and $1.25 \lambda_W$ and the distance between them is around $3 \lambda_W$, where $\lambda_W$ denotes the wavelength. Both obstacles being penetrable, the reflection of the incident wave highlighting the inclusion can be quite weak.\\

\noindent For the forward simulation, the source used $f(\xvec,t)$ is a Ricker signal located in the fluid part, with $f(\xvec,t)=(1-2\pi^2(\nu_0 t-1)^2)e^{-\pi^2(\nu_0 t-1)^2}$, $\nu_0=100 kHz$ being the central frequency corresponding to a wavelength $\lambda_W=12 mm$. Inverse problems are frequently ill-posed. Hence, an important question is the  sensitivity of the method with respect to noise in the data. Therefore, to verify the (in)sensitivity of the method with respect to the noise in the data, we added Gaussian noise by replacing the recorded data $p^S$ on $\Gamma_{SRA}$ by 
\begin{equation}
\label{labNoise}
p^S_{Noise} = (1 + \mbox{\it Coeff} * \mbox{\it randn}) * p^S
\end{equation}
where $\mbox{\it randn}$ satisfies a centered reduced normal law and the coefficient $\mbox{\it Coeff}$ is the level of noise. In our simulations,  $\mbox{\it Coeff}=10$\%. Remark that this is a multiplicative, and non an additive noise.
Indeed, following \cite{CoKr13}, this allows us to disturb the registered signal avoiding ``inverse crime'' related to
the interaction between numerical schemes used for the forward and the reverse problem. In other words,
the noise added here is a way to model the noise due to the recording device, and does not represent a background noise,
that will be rather modeled by an additive noise. In addition, also to avoid any potential inverse crime, the computational meshes used in the forward simulation are totally different from those used in the reversed ones.\\

\begin{figure}[htbp]
\centering
\includegraphics[height=0.250\textwidth]{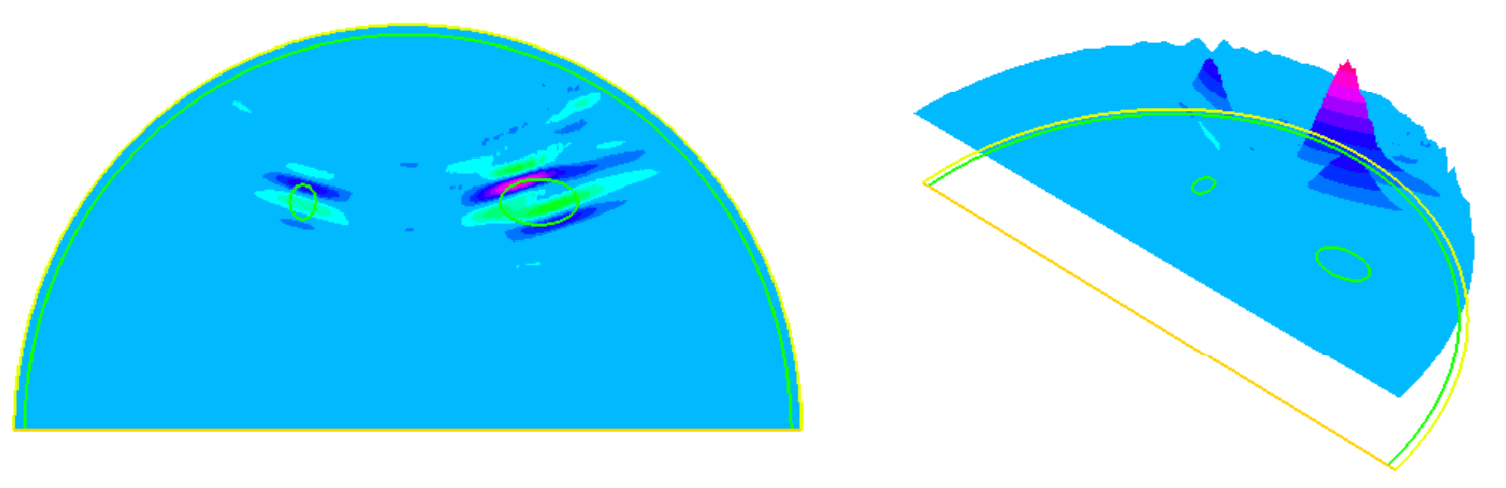}
 \caption{2D (left) and 3D (right) representation of the imaging function {\em Y-component} (\ref{RTM_equationUY}) in $\Omega_s$ for 2 malignant scatterers, with a noise level of 10\%.}
\label{Ycomp-2malign}
\end{figure}

\begin{figure}[htbp]
\centering
\includegraphics[height=0.250\textwidth]{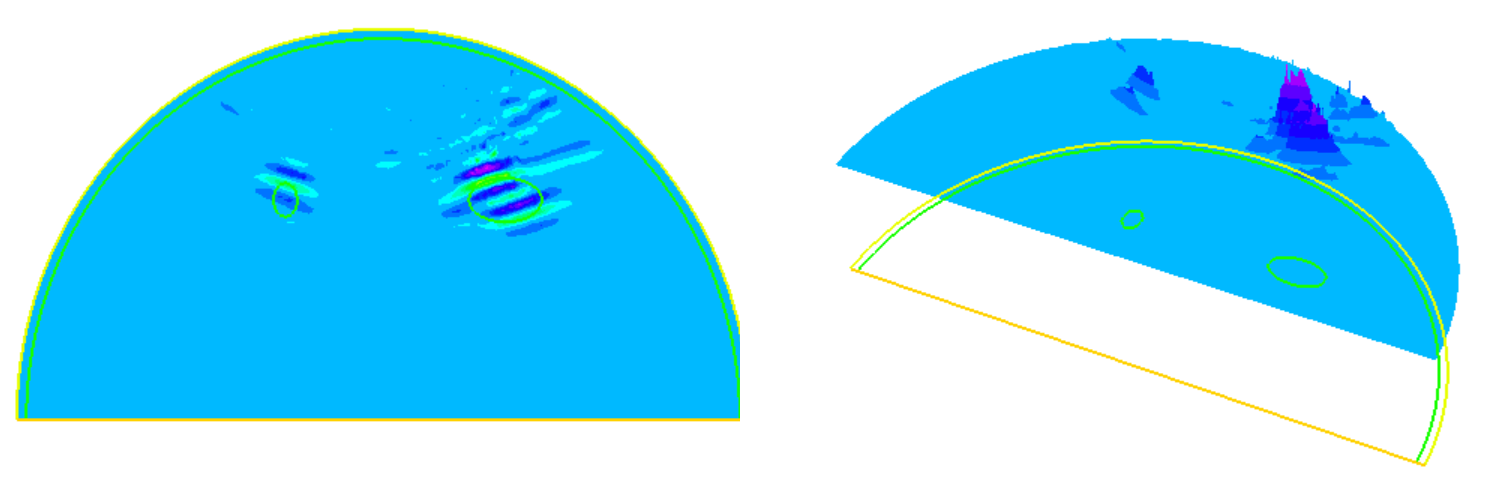}
 \caption{2D (left) and 3D (right) representation of the imaging function {\em divergence} (\ref{RTM_equationDIV}) in $\Omega_s$ for 2 malignant scatterers, with a noise level of 10\%.}
\label{Div-2malign}
\end{figure}

\noindent  In \cite{AsLi19}, we showed the feasibility of our approach in a square two-layered medium. In the following numerical tests,  we are interested to test the robustness of this method on a geometry that mimics breast tissue as described above. As a first instance, one considers the two following configurations. The first one illustrates that one can image two objects with the  same physical properties, but with different sizes, even in non homogeneous domain (fluid-solid, skin, etc.). To that end, one considers an obstacle made of two malignant elliptic scatterers $D_{1}, D_{2}$, separated by a distance of approximatively $3\lambda_W$, with major and minor axis equal to $1.25 \lambda_W, 0.75\lambda_W$ for $D_1$, and $0.58 \lambda_W, 0.4\lambda_W$ for $D_2$. In Fig.\ref{Ycomp-2malign} and Fig.\ref{Div-2malign}, we plot the value of the function RTM defined in (\ref{RTM_equation}) and (\ref{RTM_equationDIV}) respectively, in the elastic layer. The left part represents a 2D view, whereas on the right part is depicted the same image in a 3D view. This shows that one is able to determine the existence and location of the two malignant tumors.\\

%

\begin{figure}[htbp]
\begin{tabular}{lr}
{
\includegraphics[width=.4\linewidth]{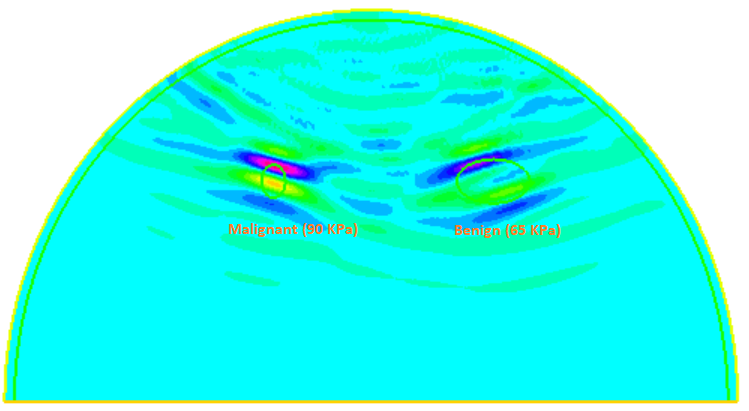}
}
&
\hspace*{2.cm}
{
\includegraphics[width=.4\linewidth]{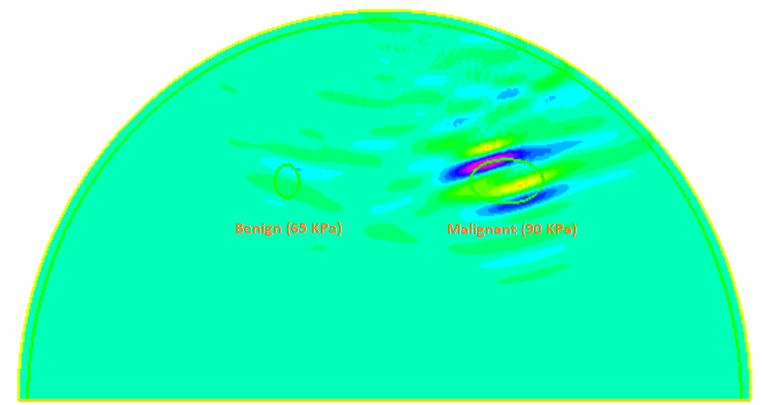}
}
\end{tabular}
 \caption{Left side: 2D representation of the imaging function {\em Y-component} (\ref{RTM_equation}) in $\Omega_s$ for 2 tumors, the left one is  malignant, the right one is benign, with a noise level of 10\%. Right side: switch between the tumors.}
\label{Ycomp-1malign-1benign}
\end{figure}

\begin{figure}[htbp]
\begin{tabular}{lr}
{
\includegraphics[width=.4\linewidth]{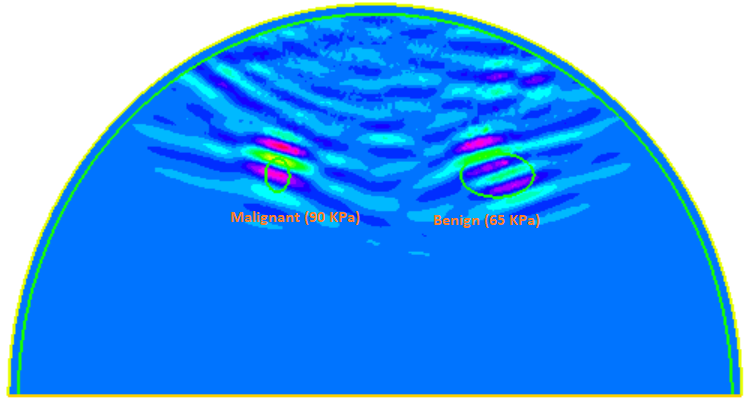}
}
&
\hspace*{2.cm}
{
\includegraphics[width=.4\linewidth]{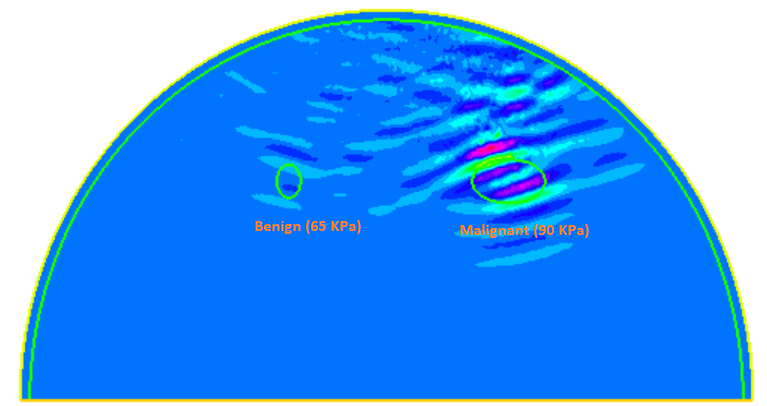}
}
\end{tabular}
 \caption{Left side: 2D representation of the imaging function {\em divergence} (\ref{RTM_equationDIV}) in $\Omega_s$ for 2 tumors, the left one is  malignant, the right one is benign, with a noise level of 10\%. Right side: switch between the tumors.}
\label{Div-1malign-1benign}
\end{figure}

\noindent The second numerical test illustrates the possibility to image two objects with  different physical properties, benign versus malignant. Here, one considers a malignant scatterer $D_1$ on the left,  and a benign one $D_2$ on the right part of $\Omega_{s}$. In addition in that case,  we reduce the size of the malignant object  to check how we can detect it. As one can see in Fig.\ref{Ycomp-1malign-1benign} and Fig.\ref{Div-1malign-1benign}, here again, the method is able to distinguish the two inclusions, and also to determine, at least qualitatively, that the one is probably malignant, whereas  the other is benign. In addition, to ensure the results are consistent, we switch between the physical properties of both tumors and confirm that the results are preserved accordingly. In the next subsection, we will introduce more quantitative criteria to determine the presence and the properties of these inclusions.\\

\noindent In both cases, the same numerical experiments have been performed with added noise on the total field, instead of the scattered one: the results are nearly the same. This is not a surprise, since it has been proved that time reversal is fairly insensitive to noise in the data (see for instance \cite{AKNT11,LTG15}). 

%

\subsection{RTM criteria to differentiate between benign and malignant tumors}

\noindent We have seen in the previous subsection that the method is able to image the inclusions, to locate them and even to (partially) ``differentiate'' them. In addition, we have also to define a mathematical criterion that tests the computed solution, to determine more quantitatively the properties of the located object. Hence, we introduce a heuristic method to assert that the tumor is benign or malignant.\\

\noindent The first one (or It) uses that fact that the Young modulus $E$,  that measures the stiffness of a solid material, is stronger in the skin than in any other tissue, benign or malignant (it is also true for the density $\rho$). Indeed, the Young modulus of the skin is between 200 KPa to 2 MPa, while it 
is between 70 KPa to 180 KPa for a malignant tumor, see for instance \cite{Fern17}. Hence, one can use it as a reference value to compute a percentage criterion that maps the values obtained on the RTM results in subsection \ref{NumResImagTR} (for criterion (\ref{RTM_equationUY}), (\ref{RTM_equationDIV}), or for another one) to a percentage that could allow us to differentiate between benign and malignant tumors. This relative RTM criterion or {\em RTM$_{percentage}$} is defined as follows:\\
For a given quantity $v$ that can be $\uvec$, $u_{2}$, $\div \uvec$, or another variable mentioned above, compute first the maximum value in space of the integral in time of $v$ for the incident wave $v^{I}$, namely
$$
\|\int_{0}^{T_f} |v^I(t, \xvec)|^{2} dt \|_{L^{\infty}(\Omega)}:= \ds\max_{\xvec \in \Omega} (\int_{0}^{T_f} |v^I(t, \xvec)|^{2} dt)
$$ 
Then, compute the normalized criterion
\begin{equation} \label{rtm_normalized}
RTM_{percentage}(\xvec) = \frac{\ds\int_{0}^{T_f} v^S_R(T_f - t, \xvec) . v^I(t, \xvec) dt\,}{\|\ds\int_{0}^{T_f} |v^I(t, \xvec)|^{2} dt \|_{L^{\infty}(\Omega)}}
\end{equation}
where $v^S_R$ denotes the time-reversed scattered of $v$.\\

\begin{figure}[htbp]
	\centering
	\begin{minipage}{.5\textwidth}
		\centering
		\includegraphics[width=.8\linewidth]{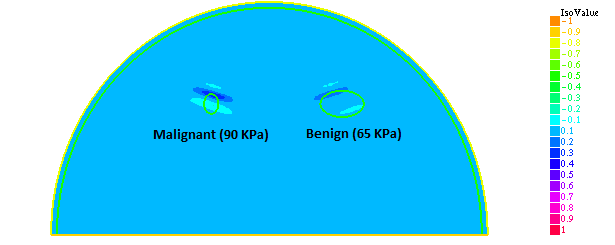}
	\end{minipage}%
	\begin{minipage}{.5\textwidth}
		\centering
		\includegraphics[width=.8\linewidth]{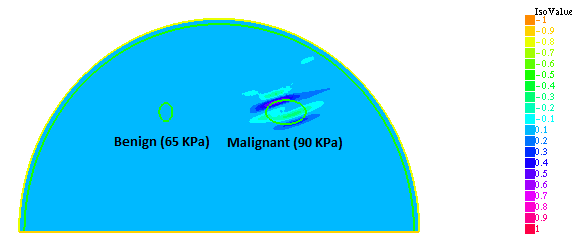}
	\end{minipage}
		\caption{RTM with normalized values (see equation \ref{rtm_normalized}) of fluid-solid medium that represents a breast tissue with 2 tumors: one benign ($E=65 KPa$) and one malignant ($E=90 KPa$); Representation of Y component component on the solid part with noise of 10\%}
	\label{fig:benign_malignat_fixed}
\end{figure}

\noindent The $RTM_{percentage}$ for $v=u_{2}$ is depicted in Fig.\ref{fig:benign_malignat_fixed}. As one can see on the right part of the figure, one can clearly identify the malignant tumor, the maximal value of the criterion being equal to $59\%$, whereas the benign one is nearly not visible. 

\noindent In addition to the previous example, this one demonstrates that TR in our context is not only able to identify the support of the scatterers, but also to discriminate (at least partially) a benign inclusion from a malignant one. However, as one can see on the left part of Fig. \ref{fig:benign_malignat_fixed} - the maximal value of the criterion being equal only to $29\%$ versus $18\%$ for the benign one - the quality of the result depends on the position of the (unique) source with respect to the scatterer to identify. 

\subsection{Work with more than a single source}

A way to improve the quality of the scatterer identification results, despite the partial information available due to the partial coverage, is to increase the number of sources to gain more coverage of the examined medium. Indeed, it has been showed in \cite{Bardos02}, \cite{BPZ02}, \cite{BPTB02}, \cite{BaPi04} or \cite{KlTi03} and \cite{CoCM00} that the spatial resolution of the refocusing is better when there is multiple scattering, that can originate from several sources. Basically, multiple scattering helps the refocusing process by causing many different waves that enrich the information contributing to the refocusing.\\

\noindent In these conditions, we have now to incorporate multiple-source measurements into the criteria introduced above.  Indeed, the RTM criteria defined in (\ref{RTM_equation}) or (\ref{rtm_normalized}) are valid when we use only one illumination. Here we propose to consider the sum of each illumination, namely to sum each criterion (\ref{rtm_normalized}) defined for a given illumination $j$. Hence, we introduce the function
\begin{equation}
\label{RTMsum}
RTM_{sum}(\xvec)=\ds \sum_{j \in SRA} \left(
 \frac{\ds\int_{0}^{T_f} v^S_R(T_f - t, \xvec,j) . v^I(t, \xvec,j) dt\,}{\|\ds\int_{0}^{T_f} |v^I(t, \xvec,j)|^{2} dt \|_{L^{\infty}(\Omega)}}
\right)
\end{equation}
where $j$ is the index of the source.\\

\noindent To illustrate our purpose,  namely to get a better coverage of the examined configuration,  we performed numerical tests by using three sources located on the SRA: the left-most element of the SRA, the middle one (used in the previous tests), and the right-most element. In each case, we consider a fluid-solid medium that represents a breast tissue with 2 tumors in the solid part (see Fig. \ref{fig:real_world_example}) with different properties.\\

\begin{figure}
	\centering
	\begin{minipage}{.5\textwidth}
		\centering
		\includegraphics[width=.8\linewidth]{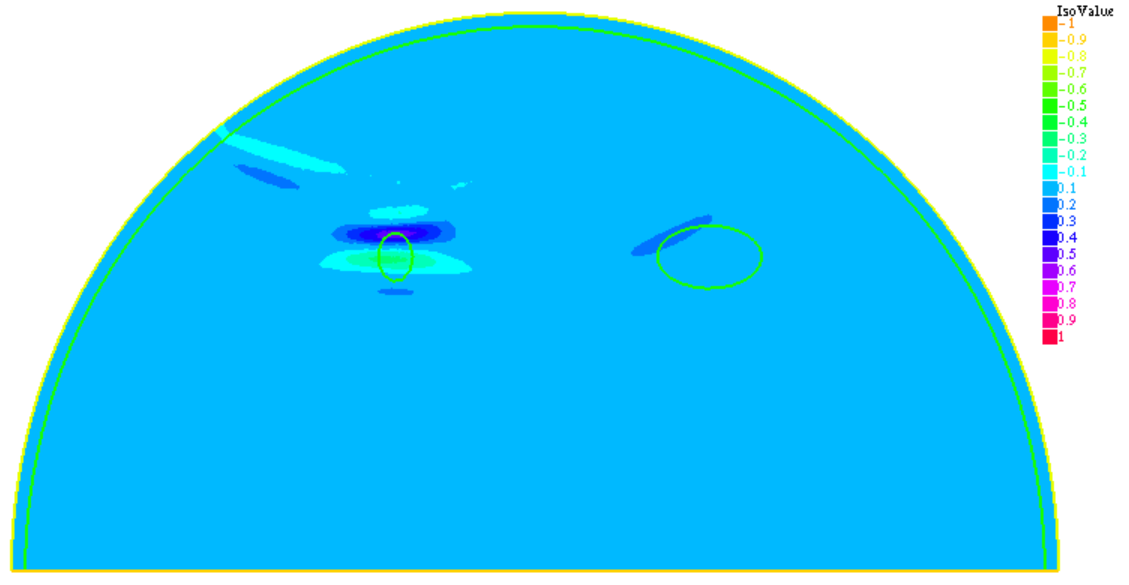}
	\end{minipage}%
	\begin{minipage}{.5\textwidth}
		\centering
		\includegraphics[width=.8\linewidth]{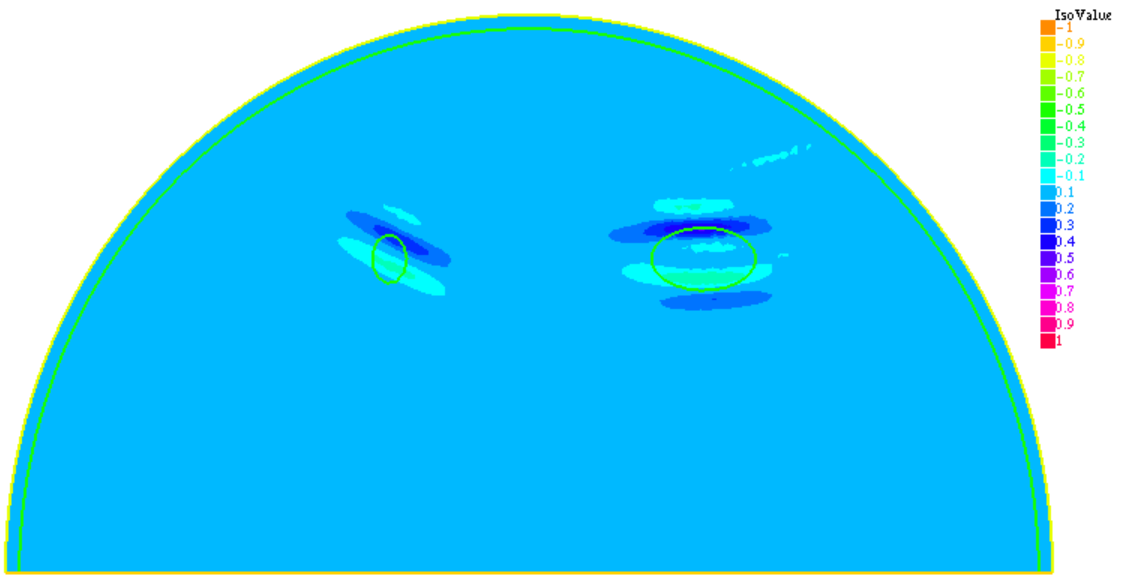}
	\end{minipage}
\caption{RTM for 2 inclusions: a small malignant tumor (90 KPa, $0.4 \lambda_{W}$) on the left, and a bigger benign one (65 Kpa, $1.25 \lambda_{W}$) on the right. Left: using the left-most source of the SRA - Right: Using the right-most source of the SRA}
\label{several-sources-smallMl}
\end{figure}

\begin{figure}
	\centering
	\begin{minipage}{.5\textwidth}
		\centering
		\includegraphics[width=.8\linewidth]{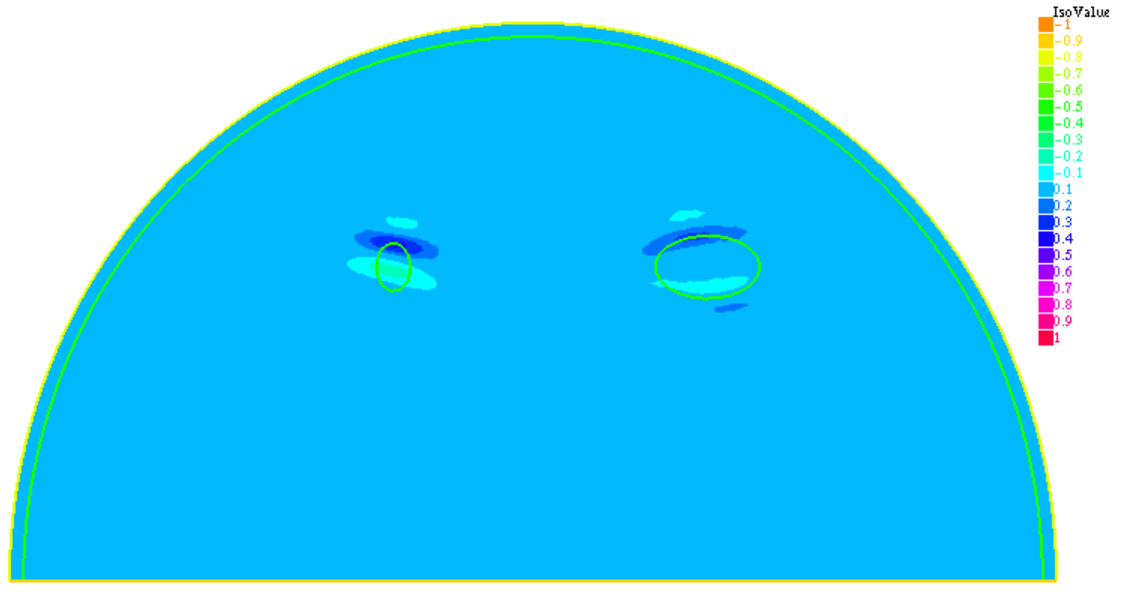}
	\end{minipage}%
	\caption{Criterion $RTM_{sum}$ (\ref{RTMsum}) obtained with 3 sources (left-most, middle, right-most) for the 2 inclusions of Figure \ref{several-sources-smallMl}.}
\label{Average-3sources-smallMl}
\end{figure}

\noindent In the first configuration (Fig.\ref{several-sources-smallMl}), a small malignant tumor (90 KPa, $0.4 \lambda_{W}$) is located on the left of the elastic domain, and a bigger benign one (65 Kpa, $1.25 \lambda_{W}$) is on the right. At the opposite, in the second configuration (Fig.\ref{Average-3sources-bigMl}), a small benign tumor (65 Kpa, $0.4 \lambda_{W}$) is on the left whereas a bigger malignant one (90 KPa, $1.25 \lambda_{W}$) is the right. In both cases, the first illumination is from the extreme left receiver, the second from the middle one (already depicted above, see Fig.\ref{fig:benign_malignat_fixed}), and the last from the extreme right receiver. The magnitude of the noise is 10\%  (see (\ref{labNoise})). Then, the result of the criterion $RTM_{sum}$ (\ref{RTMsum}) obtained with these three sources is depicted in Fig. \ref{Average-3sources-smallMl} and Fig.\ref{Average-3sources-bigMl} respectively.\\

\noindent  As expected (compare Fig. \ref{Average-3sources-smallMl} and  Fig. \ref{Average-3sources-bigMl}), it is easier to distinguish the malignant scatter from the benign one when the bigger is malignant. Indeed, the image obtained is more differentiated than in the opposite case. However, even when the benign inclusion is the bigger one,  using three sources has improved the coverage, and then improved the differentiation. In addition, it increases confidence in the result, which is less dependent on the incidence of the illumination.\\

\noindent These results are summed up in table 1, where we present the maximum value of the criterion $RTM_{sum}$ for the two configurations: for a single source (the middle one), for two sources (extreme left and right receivers) and for the three sources. It allows us to get  a quantitive criterion of discrimination between different inclusions. 

%

%
\begin{figure}
	\centering
	\begin{minipage}{.5\textwidth}
		\centering
		\includegraphics[width=.8\linewidth]{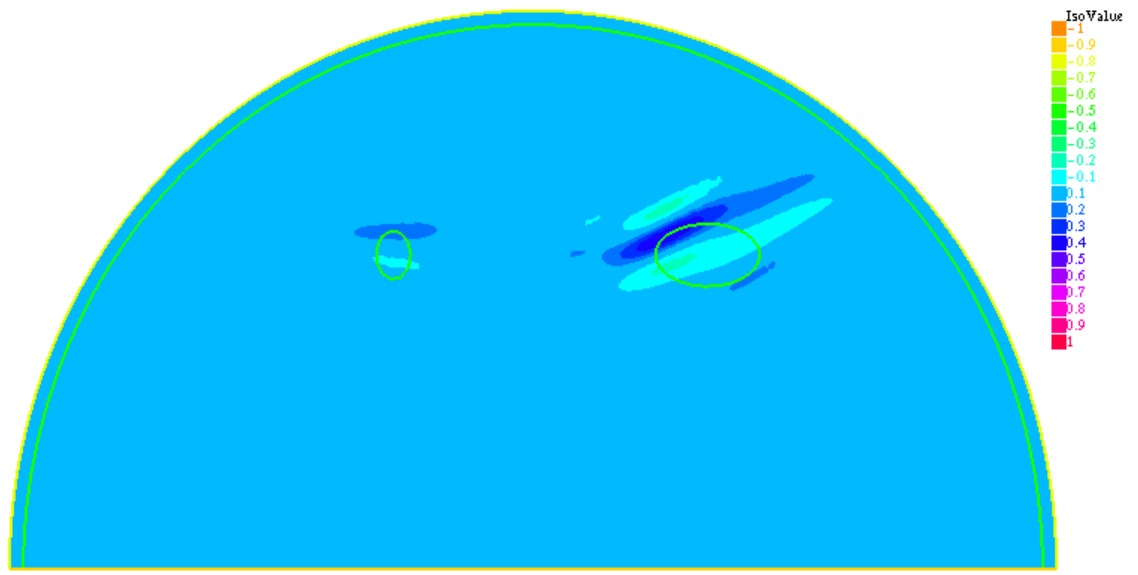}
	\end{minipage}%
	\begin{minipage}{.5\textwidth}
		\centering
		\includegraphics[width=.8\linewidth]{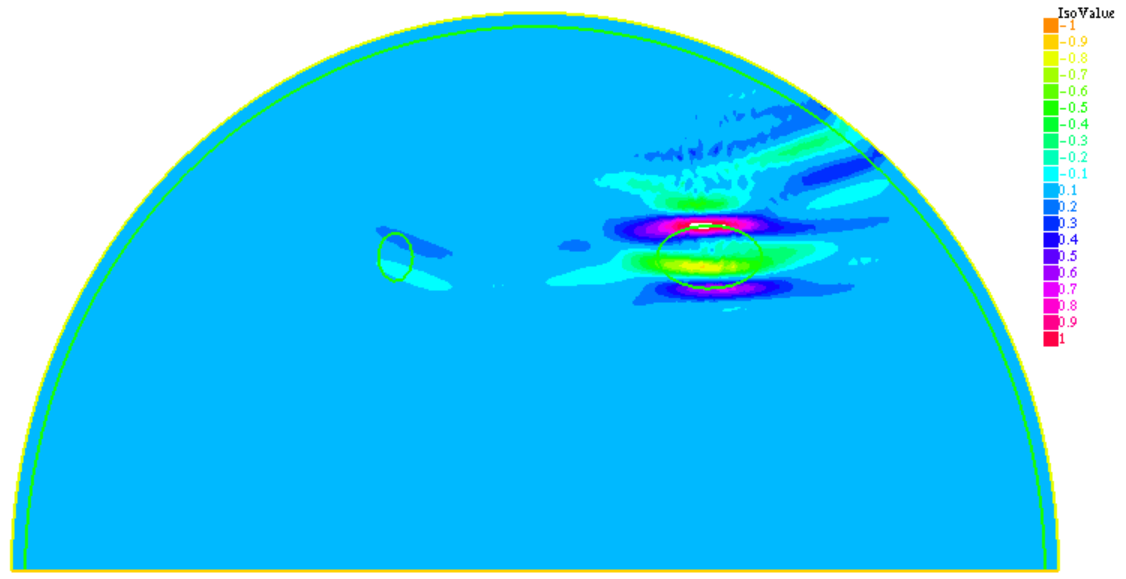}
	\end{minipage}
	\caption{RTM for 2 inclusions: a small benign tumor (65 Kpa, $0.4 \lambda_{W}$) on the left, and a bigger malignant one (90 KPa, $1.25 \lambda_{W}$) on the right. Left: using the left-most source of the SRA - Right: Using the right-most source of the SRA.}
\label{several-sources-bigMl}
\end{figure}
\begin{figure}
	\centering
	\begin{minipage}{.5\textwidth}
		\centering
		\includegraphics[width=.8\linewidth]{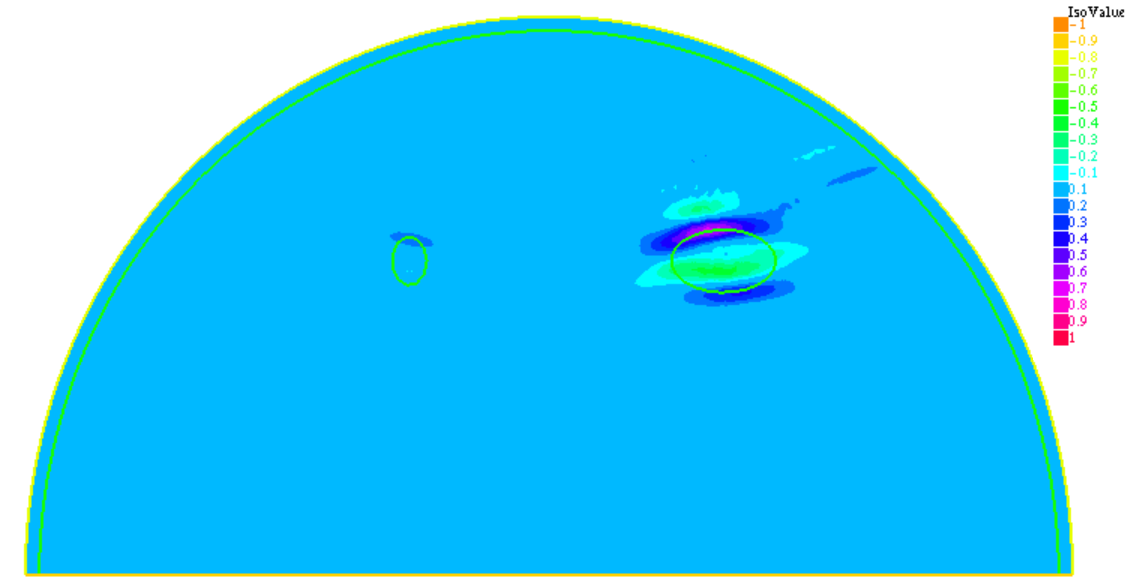}
	\end{minipage}
	\caption{Criterion $RTM_{sum}$ (\ref{RTMsum}) obtained with 3 sources (left-most, middle, right-most) for the 2 inclusions of Figure \ref{several-sources-bigMl}.}
\label{Average-3sources-bigMl}
\end{figure}

\begin{table}[!h]
\label{table1}
\centering
\begin{tabular}{|c|c|c|c|}
\hline
	$\|RTM_{sum}(\xvec)\|_{L^\infty(\Omega)}$ & 1 source & 2 sources & 3 sources \\
	\hline
	small Ml  \# big Bn  & 0.291 \# 0.181  & 0.384 \# 0.286	&  0.351 \# 0.262 \\
	\hline
	small Bn \# big Ml  &0.109 \# 0.595 & 0.19 \# 0.782
	&  0.192  \# 0.706  \\
	\hline
\end{tabular}
\caption{Maximum value of the criterion $RTM_{sum}$ defined on (\ref{RTMsum}) for a fluid-solid medium that represents a breast tissue with 2 tumors with different properties (noise of 10\%). In the first row, the small tumor is malignant (MI) with $E=90 KPa$ and the big one is benign (Bn) with $E=65 KPa$. In the second row, the small tumor is benign whereas the big one is malignant.}
\end{table}

\newpage
\subsection{Improving the coverage}

\begin{figure}[h]
	\centering
	\includegraphics[width=.4\linewidth]{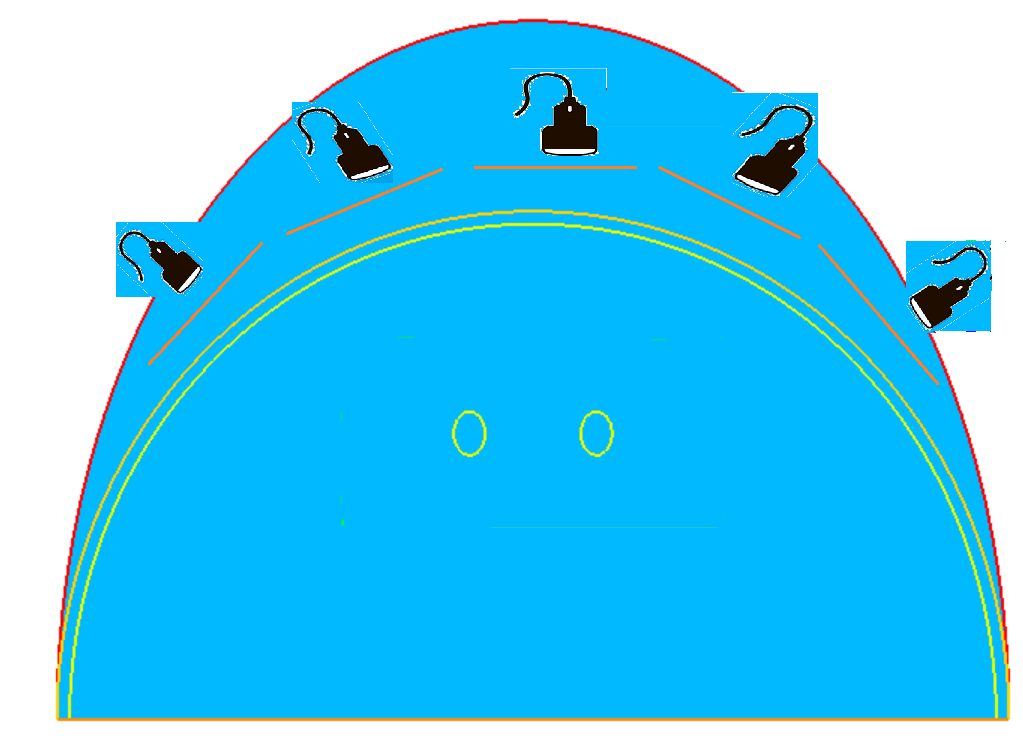}
	\caption{Multiple SRA that cover the medium, each SRA consist of 27 sensors and is sequentially placed in positions from 1 to 5.}
	\label{Fig-several-probes}
\end{figure}

In this last part, we propose to study the improvement of the coverage by using several sources with different angles, with the hope of enhancing the imaging results. To that end, we will use several (narrow) SRAs that will produce images from different  angles, and then combine them. This process is actually similar to experimental ultrasound examination using several probes that has partial coverage of the examined tissue, see Fig. \ref{Fig-several-probes}.\\

\begin{figure}[h]
	\centering
	\begin{minipage}{.49\textwidth}
		\centering
	\includegraphics[width=.8\linewidth]{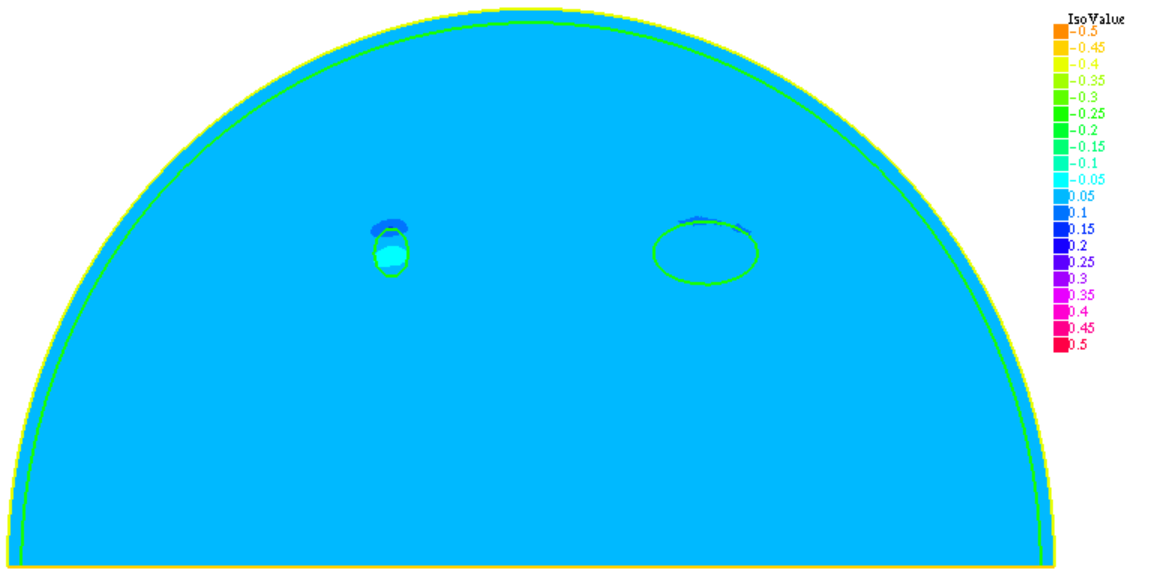}
	\end{minipage}
	\begin{minipage}{.49\textwidth}
		\centering
		\includegraphics[width=.8\linewidth]{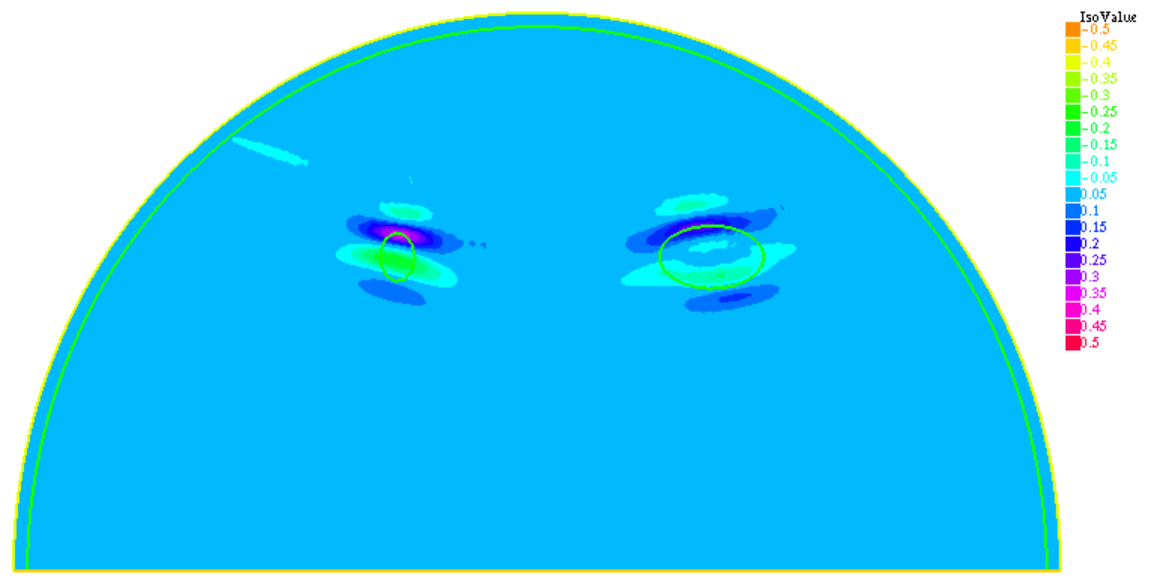}
	\end{minipage}
	\caption{Imaging function $RTM_{sum}$ of the Y-component. The left scatter is malignant (90 Kpa), the right one is benign (65 KPa). Left image:  with 5 SRAs - Right image: same as Fig. \ref{Average-3sources-smallMl}.}
\label{5-probes-vs-3SRA-malign-benign}
\end{figure}

\noindent Hence, we perform the same experiment as above, namely with a small malignant tumor (90 KPa, $0.4 \lambda_{W}$)  on the left, and a bigger benign one (65 Kpa, $1.25 \lambda_{W}$) on the right (see Fig.\ref{Average-3sources-smallMl}), but this time with 5 ``probes'' as described above. The imaging function is represented in Fig. \ref{5-probes-vs-3SRA-malign-benign} (left part), where we reproduced (right part) the previous case (with the same scale) to facilitate the comparison. As one can see, using multiple measurements with several incidences improve the differentiation between the 2 scatterers, the right benign one becoming almost invisible in the left part of the figure, whereas on the right image depicted with the same scale, it remains quite difficult to differentiate the 2 scatterers. This example shows that with more illuminations and more ``incidences'', one can improve the differentiation between the malignant and the benign, even in the ``difficult'' case of a small malignant tumor versus a bigger benign one.\\

\begin{figure}[h]
	\centering
	\begin{minipage}{.49\textwidth}
		\centering
	\includegraphics[width=.8\linewidth]{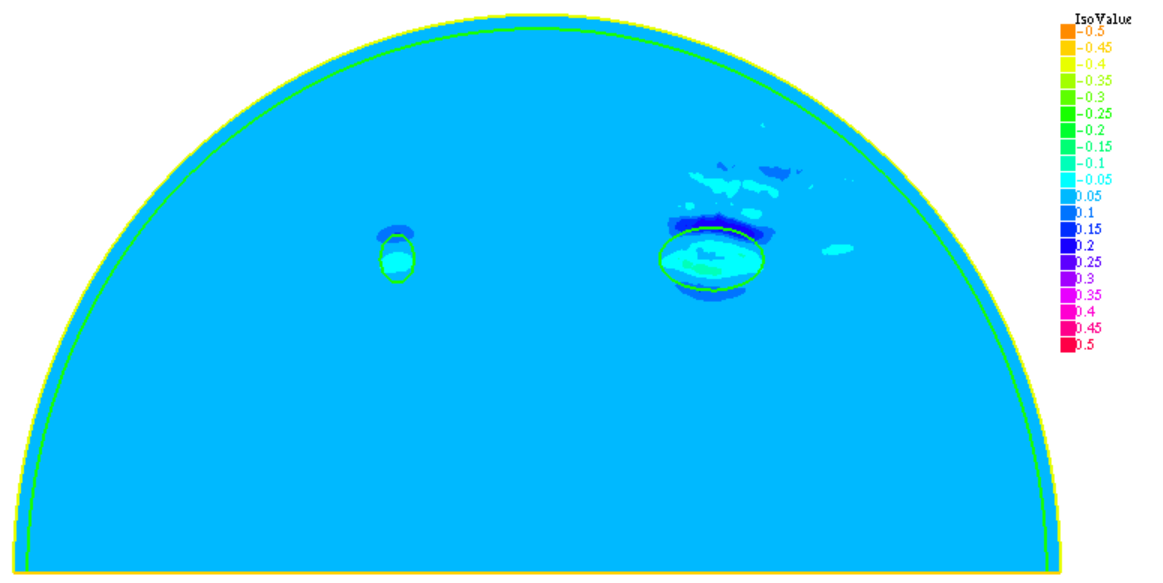}
	\end{minipage}
	\begin{minipage}{.49\textwidth}
		\centering
	\includegraphics[width=.8\linewidth]{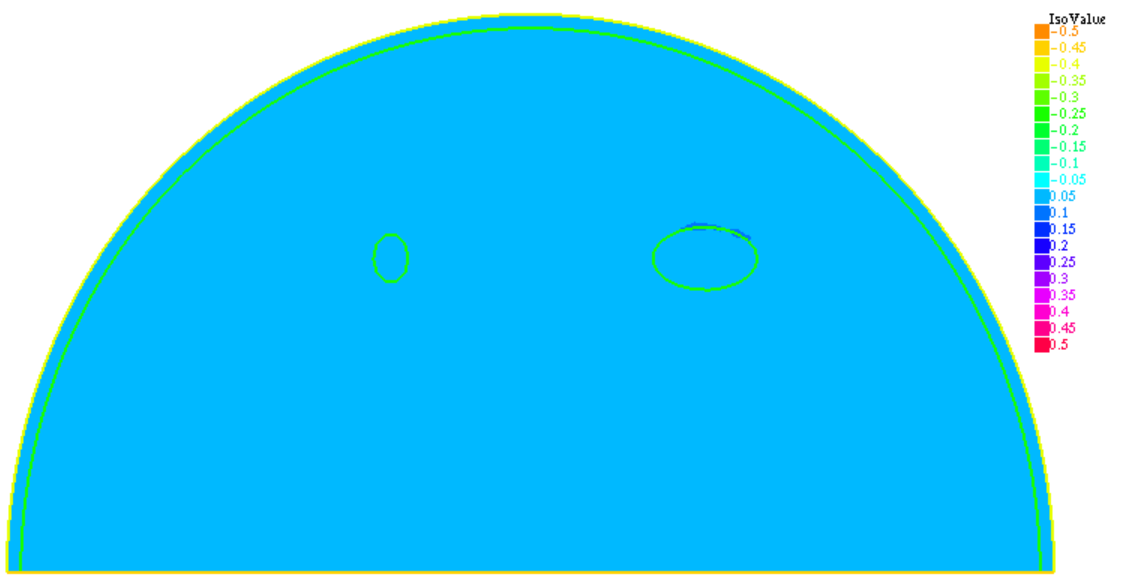}	
	\end{minipage}
	\caption{Imaging function $RTM_{sum}$ of the Y-component when both tumors are malignant (90 KPa) (left figure) compared to the case where both tumors are benign (65 KPa) (right figure).}
\label{5-probes-malign-benign}
\end{figure}

\noindent Now, in this second numerical test, we want to investigate the case where both scatterers have the same physical properties.  We consider the same geometry as above, but in that case, both scatterers are either malignant or benign objects. Here again, we use, as before, 5 illuminations from 5 different incidences. Results are depicted in Fig. \ref{5-probes-malign-benign} with the same scale. This shows that the malignant objects are clearly visible, even both are malignant, whereas, with the same scale, the benign ones are almost invisible. This illustrates the possibility to retrieve the properties of the tumors not only in a ``relative'' way, but also in a more ``absolute'' way.

\begin{figure}[h]
	\centering
	\includegraphics[width=.4\linewidth]{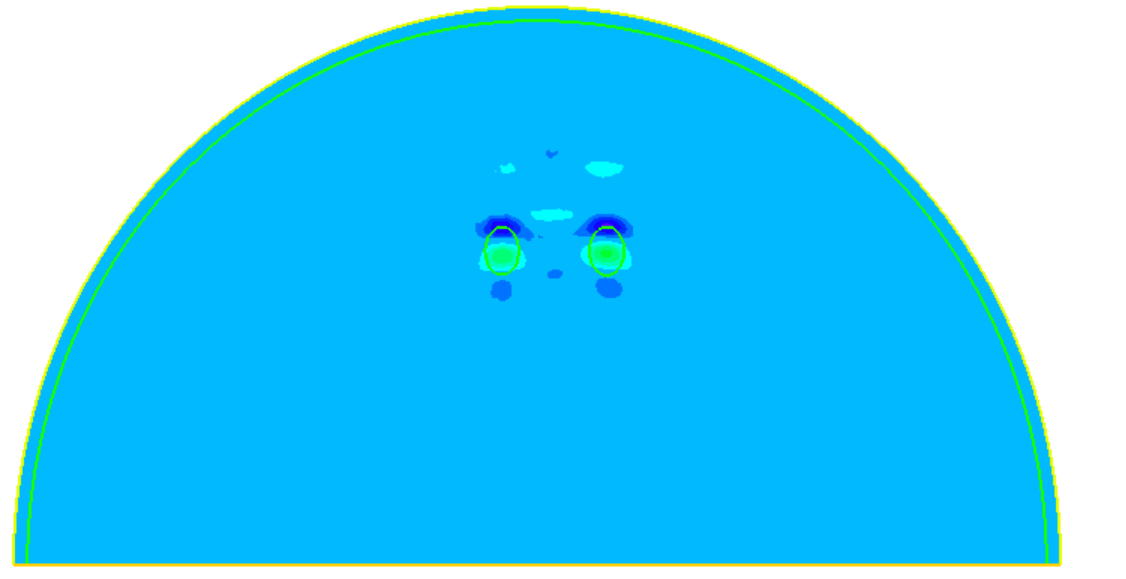}
	\caption{Imaging function $RTM_{sum}$ of the Y-component for 2 close malignant (90 KPa) objects.}
\label{5-probes-malign-malign-close}
\end{figure}

\noindent As a last illustration, we investigate the ability of the method to differentiate between two close inclusions. In this numerical example, we consider 2 malignant scatterers with the same size of about $0.4 \lambda_W$, separated by a distance of $\frac{\lambda_W}{2}$. Figure \ref{5-probes-malign-malign-close} displays results of the imaging function $RTM_{sum}$ for the Y-component. One can clearly observe the 2 inclusions in the medium. Hence, even it the tumors are quite close one to the other, the method is able to image them and even to differentiate them.

\section*{Conclusion} 
\label{Conclusion}

\noindent In this paper, we proposed a time-reversal method for scatterer identification in a non homogeneous elastic medium. Because of the concerned configuration, i.e. a layered fluid-solid medium, one is forced to consider only partial aperture case with data  recorded only in the acoustic part of the medium.\\

\noindent  Since the medium is not homogeneous (fluid layer, skin, elastic medium,..), and due to the presence of noise, the time-reversal method is not able to ``recreate the past'' of the wave propagation, even if the location and the properties of the inclusions which diffracted the signals were know: this due to the non-homogeneity of the medium.\\

\noindent  However, we shown that these data contain enough information to identify the scatterers, that is the properties and locations of the inclusions backscattering the signal. This was made possible by introducing criteria, derived from the reverse time migration framework, to construct images of the inclusions and to determine their locations.\\

\noindent The dependence of the approach to several parameters (aperture, number of sources, etc.) was also investigated. In particular, it was shown that one can differentiate between a benign and malignant close inclusions, that are treated following the same procedure. We also shown that the method is able to differentiate between one and two inclusions, even quite close. In addition, this approach has proved to be fairly insensitive to noise in the data, as usual in time-reversal methods. Finally, although our numerical examples are in two-space dimensions for fluid-solid media, they can be extended to 3D or other wave propagation equations, especially for linear elastodynamics, without fluid part. In the same spirit,  possible extensions to the harmonic case (see example in \cite{AKNT11} for the acoustic case) seems to be within reach.


%
\end{document}